\documentclass{article}
\usepackage{amssymb}
\usepackage{amsmath}
\usepackage{amsfonts}

\newtheorem{theorem}{Theorem}

\newtheorem{definition}[theorem]{Definition}

\newtheorem{proposition}[theorem]{Proposition}

\input{tcilatex}
\begin{document}

\title{Curvature without connection}
\author{Erc\"{u}ment H. Orta\c{c}gil}
\maketitle

\begin{abstract}
We show that the alternative theory of Lie groups and geometric structures
proposed in the recent book [Or1] can be developed independently of
connections. We show the details of this connection-free approach in the
cases of absolute parallelism, affine and Riemannian structures and outline
the method in the general case.
\end{abstract}

\section{Introduction}

A nonlinear/linear prehomogeneous geometry (PHG) as defined in [Or1] is a
special transitive Lie equation in finite and infinitesimal forms whose
theory is developed by D.C.Spencer and his coworkers around 1970. Lie
equations, whose origin traces back to the original works of Sophus Lie, are
later generalized to groupoids (nonlinear or finite in the old terminology)
and algebroids (linear or infinitesimal). In [Or1] we defined the curvature
of a PHG under a technical assumption which is equivalent to this curvature
being the curvature of a connection. This turns out to be a very strong
condition as explained in [Or1] (see pages 188-192). As a solution, we
proposed another definition based on an idea communicated to us by A. Blaom.
With time, however, we came to realize that this second definition that we
advocated as the ultimate solution throughout [Or1] has the same weakness as
the classical definition of a Cartan geometry: It is difficult to construct
nonflat examples (see [CS] for parabolic geometries). The purpose of this
paper is to give a detailed analysis of an affine and Riemannian PHG showing
that the above mentioned technical assumption is redundant, a fact stated
also in [Or1] without proof. This paper can be regarded also as a revision
and extension of [Or2].

A few words about the general philosophy of [Or1] and this paper are in
order here: Principal bundles, vector bundles and connections on these
bundles are essentially topological objects. It is a standard practice today
to do geometry using these topological tools. By geometry we mean here the
study of transitive actions of Lie groups in the spirit of Klein's Erlangen
Program complemented with a concept of curvature, or briefly, the study of
"symmetry deformed by curvature" in the words of Blaom ([Bl]). In the
approach commonly accepted today, which regards connection as the primary
object, the search for some "special connection" that suits the geometric
structure "best" becomes an unavoidable and sometimes an ardous task ([CS]).
Naturally the question arises whether it is possible to start directly with
geometry building on some simple and intuitive principles based on "group
action" and recover the above topological framework as a generalization.
Such a geometric framework is proposed in [Or1]. In this alternative
approach certain "splittings" replace "special connections". These
splittings are incorporated into the definition of the geometric structures
and they do not always define connections, but when they do, these
connections turn out to be "very special".

This alternative approach has many advantages over the classical one. For
instance,

1) As mentioned above, it frees us from the search for special connections,
and consequently from one of the most intriguing concepts in differential
geometry: torsion.

2) It is much easier to construct nonflat examples.

3) It is modelfree, that is, the flat model (which is some homogeneous
space) is not fixed beforehand. In particular, it avoids the Maurer-Cartan
form which fixes the model. This fact gives an immense conceptual depth to
the theory. For instance, Poincare Conjecture becomes tantamount to proving
the existence of a flat absolute parallelism on a compact and simply
connected 3-manifold.

4) A PHG is by definition a first order mixed system of PDE and curvature is
by definition the integrability conditions of this system. Therefore, this
definition of curvature is both elementary and intuitive and avoids the
subtle aspects of Spencer cohomology, like involution, acyclicity, formal
integrability...etc.

5) All geometric structures are defined on equal logical footing.

6) Nevertheless, there is an hiearchy determined by the order of jets
measuring the inner complexity of the structure. In this hierarchy, the
simplest one is absolute parallelism which renders in the flat case the
theory of \textit{simply transitive }actions of Lie groups and Lie algebras.

7) The study of \textit{all }PHG's reduces to the case of parallelism and
the principle of this reduction is the moving frame method introduced by
Fels and Olver in [FO1], [FO2] as a correction and perfection of Cartan's 
\textit{repere mobile. }This method equips the present approach with
powerful computational algorithms.

It is our hope that the differential geometers will find this alternative
approach worthwhile.

\section{Review of absolute parallelism}

We started Part 1 of [Or1] with a quote from Einstein: \textit{Everything
should be made as simple as possible, but not simpler.} With time, however,
we came to realize that Part 1 of [Or1] is in need of further
simplification, especially the unnecessarily long and obscure proof of the
crucial Proposition 6.8. This simplification has a striking consequence: The
connections $\nabla ,\widetilde{\nabla },$ which seem to play a fundamental
role in Parts 1, 2 of [Or1], emerge as consequences of more fundamental
concepts and \textit{one can develop the whole theory without even
mentioning them!} \textit{\ More importantly, this new approach applies word
by word to all geometric structures. }The key fact is that the theory in
Part 1 of [Or1] can be based on the "structure object" which is passed over
in silence in [Or1] except in Chapter 7. Our purpose here is to recast Part
1 of [Or1] in this simplified form, refering to [Or1] for some further
technical details.

Let $M$ be a smooth manifold and $\varepsilon $ a splitting of the groupoid
projection $\pi :\mathcal{U}_{1}\rightarrow \mathcal{U}_{0}=M\times M.$ It
is easy to show that such $\varepsilon $ exists if and only if $M$ is
parallelizable ([Or1], Proposition 1.2) . Since $\varepsilon $ is a
homomorphism of groupoids, it satisfies the following identities:

\begin{equation}
\varepsilon _{a}^{i}(y,z)\varepsilon _{j}^{a}(x,y)=\varepsilon _{j}^{i}(x,z)
\end{equation}%
\begin{equation}
\varepsilon _{j}^{i}(x,x)=\delta _{j}^{i}
\end{equation}%
\begin{equation}
\varepsilon _{a}^{i}(y,x)\varepsilon _{j}^{a}(x,y)=\varepsilon
_{j}^{i}(x,x)=\delta _{j}^{i}
\end{equation}%
for all $x,y,z\in M.$ We fix a base point $e\in M$ and define a geometric
object $w$ with components $(w_{j}^{i}(x))$ on $(U,x)$ by

\begin{equation}
w_{j}^{i}(x)\overset{def}{=}\varepsilon _{j}^{i}(e,x)
\end{equation}

Now (1), (3) and (4) give

\begin{equation}
\varepsilon _{j}^{i}(x,y)=\varepsilon _{a}^{i}(e,y)\varepsilon
_{j}^{a}(x,e)=w_{a}^{i}(y)\widetilde{w}_{j}^{a}(x)\text{ \ \ \ \ \ \ }%
\widetilde{w}=w^{-1}
\end{equation}

Since $\varepsilon (x,y)\circ \varepsilon (e,x)=\varepsilon (e,y)$ by (1),
the components $(w_{j}^{i}(x))$ are subject to the transformation law 
\begin{equation}
\frac{\partial y^{i}}{\partial x^{a}}w_{j}^{a}(x)=w_{j}^{i}(y)
\end{equation}%
upon a coordinate change $(x)\rightarrow (y)$ in view of (4). Therefore, if
we are given an object $w$ on $M$ whose components $(w_{j}^{i}(x))$
transform according to (6), we can define $\varepsilon $ by (5) which is
easily seen to be a splitting and therefore satisfies (1)-(3). Conversely,
given a splitting $\varepsilon $ satisfying (1)-(3), we can fix a base point 
$e$ and define $w$ by (4) whose components transform according to (5). We
observe that $w$ determines $\varepsilon $ canonically whereas $\varepsilon $
determines $w$ modulo the choice of a base point $e$ and some coordinates
around $e$ satisfying $w_{j}^{i}(e)=\delta _{j}^{i}.$ It follows that $w$ is
more intrinsic than $\varepsilon $ and henceforth, we will concentrate on $w$
rather than $\varepsilon .$

\begin{definition}
The geometric object $w=(w_{j}^{i}(x))$ with the transformation law (6) is
the structure object of the parallelizable manifold $(M,w).$
\end{definition}

There is a more conceptual derivation of $w:$ Let $G$ be an abstract Lie
group and $L_{g}:G\rightarrow G$ be the left translation $L_{g}(x)=gx.$ Let $%
\mathcal{U}_{1}^{e,e}$ the set of $1$-jets of \textit{all }local
diffeomorphisms (which we call $1$-arrows) that fix $e.$ Then $\mathcal{U}%
_{1}^{e,e}$ is a Lie group and a choice of coordinates around $e$ identifies 
$\mathcal{U}_{1}^{e,e}$ with $GL(n,\mathbb{R)}.$ The only $L_{g}$ that fixes 
$e$ is $L_{e}$ which we identify with its $1$-jet $Id$ in $\mathcal{U}%
_{1}^{e,e}\cong GL(n,\mathbb{R)}.$ Now $GL(n,\mathbb{R)}$ acts on the left
coset space $GL(n,\mathbb{R)}/\{Id\}=GL(n,\mathbb{R)}$ by matrix
multiplications as $(w_{j}^{i})\rightarrow g_{a}^{i}w_{j}^{a}$ and this
group action can be used to define the structure object $w.$ This derivation
applies to all homogeneous spaces $G/H$ and can be used to define the
"structure object of that particular geometry" as we will see below.

We now define%
\begin{equation}
\Gamma _{jk}^{i}(x)\overset{def}{=}\left[ \frac{\partial \varepsilon
_{k}^{i}(x,y)}{\partial y^{j}}\right] _{y=x}=\frac{\partial \varepsilon
_{a}^{i}(y,x)}{\partial x^{j}}\varepsilon _{k}^{a}(x,y)=\frac{\partial
w_{a}^{i}(x)}{\partial x^{j}}\widetilde{w}_{k}^{a}(x)
\end{equation}

The second equality in (7) follows by differentiating (1) with respect to $z$
at $z=x$ and the third equality follows from (5). In particular, note that
the third expression in (7) is independent of $y.$ Similarly, we obtain

\begin{equation}
-\Gamma _{jk}^{i}(x)=\left[ \frac{\partial \varepsilon _{k}^{i}(x,y)}{%
\partial x^{j}}\right] _{y=x}=\varepsilon _{a}^{i}(y,x)\frac{\partial
\varepsilon _{k}^{a}(x,y)}{\partial x^{j}}=w_{a}^{i}(x)\frac{\partial 
\widetilde{w}_{k}^{a}(x)}{\partial x^{j}}
\end{equation}

Using (6), we now define a subgroupoid $\mathcal{N}_{1}(w)\subset \mathcal{U}%
_{1}$ as follows: The fiber $\mathcal{N}_{1}(w)^{\overline{x},\overline{y}}$
of $\mathcal{N}_{1}(w)$ over $(\overline{x},\overline{y})\in \mathcal{U}%
_{0}=M\times M$ consists of those $1$-arrows $(\overline{x},\overline{y},%
\overline{f}_{1})=(\overline{x}^{i},\overline{y}^{i},\overline{f}_{j}^{i})$
of $\mathcal{U}_{1}$ from $\overline{x}$ to $\overline{y}$ that preserve the
geometric object $w,$ that is

\begin{equation}
\overline{f}_{a}^{i}w_{j}^{a}(\overline{x})=w_{j}^{i}(\overline{y})
\end{equation}

\begin{definition}
The subgroupoid $\mathcal{N}_{1}(w)\subset \mathcal{U}_{1}$ is the
invariance groupoid of the structure object $w.$
\end{definition}

$\mathcal{N}_{1}(w)$ is the most important example of a prehomogeneous
geometry (PHG) and historically it is called an absolute parallelism. We
observe that if $(\overline{x},\overline{y},\overline{f}_{1})=(\overline{x}%
^{i},\overline{y}^{i},\overline{f}_{j}^{i})$ is in $\mathcal{N}_{1}(w),$
then $\overline{f}_{j}^{i}=$ $w_{a}^{i}(\overline{y})\widetilde{w}_{j}^{a}(%
\overline{x})$ by (9) and therefore $\overline{f}_{j}^{i}$ is determined by $%
\overline{x}^{i}$ and $\overline{y}^{i}.$ Therefore, above any $0$-arrow $(%
\overline{x},\overline{y})\in \mathcal{U}_{0}=M\times M,$ there is a unique $%
1$-arrow $(\overline{x},\overline{y},\overline{f}_{1})$ of $\mathcal{N}%
_{1}(w)$ which we write as $\varepsilon (\overline{x},\overline{y})=(%
\overline{x},\overline{y},\overline{f}_{1}).$ We define $\mathcal{N}_{0}(w)%
\overset{def}{=}\mathcal{U}_{0}=$ the pair groupoid $M\times M$ so that $%
\mathcal{N}_{0}(w)\cong \mathcal{N}_{1}(w)$ where the isomorphism $\cong $
of groupoids is given by (omitting the indices) $(\overline{x},\overline{y}%
)\rightarrow (\overline{x},\overline{y},w(\overline{y})\widetilde{w}(%
\overline{x}))=(\overline{x},\overline{y},\varepsilon (\overline{x},%
\overline{y})).$ Therefore $\mathcal{N}_{1}(w)=\varepsilon (\mathcal{U}%
_{0})\subset \mathcal{U}_{1}.$

\begin{definition}
A local bisection of $\mathcal{U}_{1}$ consists of

1) An arbitrary local diffeomorphism $f:U\rightarrow V=f(U)$

2) A smooth choice of $1$-arrows of $\mathcal{U}_{1}$ from $x\in U$ to $%
f(x)\in V$
\end{definition}

In coordinates, a local bisection of $\mathcal{U}_{1}$ is of the form

\begin{equation}
(x^{i},f^{i}(x),f_{j}^{i}(x))
\end{equation}

If we fix $x=\overline{x}$ in (10), then (10) becomes the $1$-arrow denoted
by $(\overline{x},\overline{y},\overline{f}_{1})$ above. Similarly we define
a local bisection of $\mathcal{N}_{1}(w)$ by requiring the $1$-arrows in
Definition 3 to belong to $\mathcal{N}_{1}(w).$ In this case, note that
there is no "choice" of $1$-arrows since from $x$ to $f(x)$ there is a
unique $1$-arrow of $\mathcal{N}_{1}(w).$

If (10) is a local bisection of $\mathcal{N}_{1}(w),$ then it satisfies

\begin{equation}
f_{a}^{i}(x)w_{j}^{a}(x)=w_{j}^{i}(f(x))
\end{equation}
according to (9).

\begin{definition}
The local bisection (10) is prolonged (or holonomic) if
\end{definition}

\begin{equation}
f_{j}^{i}(x)=\frac{\partial f^{i}(x)}{\partial x^{j}}
\end{equation}

Now $\mathcal{N}_{1}(w)$ defines a first order nonlinear system of PDE's on
the pseudogroup $Diff_{l}(M)$ of local diffeomorphisms of $M$ as follows:
Some local diffeomorphism $y=f(x)$ is a local solution of $\mathcal{N}%
_{1}(w) $ if it satisfies

\begin{equation}
\frac{\partial f^{i}(x)}{\partial x^{a}}w_{j}^{a}(x)=w_{j}^{i}(f(x))
\end{equation}%
or equivalently

\begin{equation}
\frac{\partial f^{i}(x)}{\partial x^{j}}=\varepsilon _{j}^{i}(x,f(x))
\end{equation}

Therefore, a local solution of $\mathcal{N}_{1}(w)$ is a prolonged bisection
of $\mathcal{N}_{1}(w)$ and now the question is whether $\mathcal{N}_{1}(w),$
which clearly admits many local bisections (as we can choose $y=f(x)$
arbitrarily), admits any prolonged bisections. The fibers $\mathcal{N}%
_{1}(w)^{\overline{x},\overline{y}}$ of $\mathcal{N}_{1}(w)$ now serve as
the initial conditions for the PDE (13): For any initial condition $(%
\overline{x},\overline{y},\overline{f}_{1})$ satisfying (9), can we find a
prolonged bisection $(x,f(x),\frac{\partial f}{\partial x})$ defined around $%
\overline{x}$ and satisfying $f(\overline{x})=\overline{y},$ $\frac{\partial
f}{\partial x}(\overline{x})=\overline{f}_{1}$ $?$ This is a particular
instance of a general framework of PDE's: Jets (1-arrows in our case) are
pointwise but not always locally derivatives. A PDE (for instance (11)) is a
condition on jets. Solving the PDE is equivalent to replacing jets by
derivatives (or replacing (11) by (13))

Now the well known existence and uniqueness theorem for first order systems
of PDE's with initial conditions asserts that the answer to the above
question is affirmative if and only if the integrability conditions of (13)
are satisfied. The check these integrability conditions, we differentiate
(13) with respect to $x^{k},$ substitute back from (13) and alternate $j,k.$
After some computation, this condition turns out to be

\begin{equation}
\left[ \frac{\partial w_{a}^{i}(f(x))}{\partial x^{j}}\widetilde{w}%
_{k}^{a}(f(x))\right] _{[jk]}-f_{d}^{i}(x)\left[ \frac{\partial w_{c}^{d}(x)%
}{\partial x^{a}}\widetilde{w}_{b}^{c}(x)\right]
_{[ab]}g_{j}^{a}(x)g(x)_{k}^{b}=0
\end{equation}%
for \textit{all }bisections $(x,y,f_{1})=(x^{i},f^{i}(x),f_{j}^{i}(x))$ of $%
\mathcal{N}_{1}(w)$ where $g_{1}=f_{1}^{-1}.$ We now define the geometric
object $I(w)$ by defining its components $I_{jk}^{i}(w;x)$ on $(U,x)$ by

\begin{equation}
I_{jk}^{i}(w;x)\overset{def}{=}\left[ \frac{\partial w_{a}^{i}(x)}{\partial
x^{j}}\widetilde{w}_{k}^{a}(x)\right] _{[jk]}=\left[ \Gamma _{jk}^{i}(x)%
\right] _{[jk]}
\end{equation}

\begin{definition}
$I(w)$ is the integrability object of $w.$
\end{definition}

The name for $I(w)$ will be justified below. In [Or1], $I(w)$ is called the
torsion of $\mathcal{N}_{1}(w)$ because it turns out to be the torsion of a
linear connection on the tangent bundle $T(M)\rightarrow M.$ However, this
is a very misleading terminology as it holds only in the case of absolute
parallelism as we will see below.

Using the LHS of (15) and (16), we define $\mathcal{R}_{jk}^{i}(x,y)$ on $%
M\times M$ by

\begin{equation}
\mathcal{R}_{jk}^{i}(x,y)\overset{def}{=}%
I_{jk}^{i}(w;y)-I_{ab}^{c}(w;x)f_{c}^{i}g_{j}^{a}g_{k}^{b}
\end{equation}%
for any initial condition $(x,y,f_{1})\in \mathcal{N}_{1}(w)^{x,y},$ $%
g_{1}=f_{1}^{-1}.$

\begin{definition}
$\mathcal{R}_{jk}^{i}(x,y)$ is the nonlinear curvature of $\mathcal{N}%
_{1}(w).$
\end{definition}

$\mathcal{R}$ is denoted by $\overleftarrow{\mathcal{R}}$ in [Or1] (2.12,
pg.20). Note that $f_{1}=g_{1}^{-1}$ on the RHS of (17) is not an
independent variable since it is determined by $x,y.$ Now if (13) admits
local solutions with arbitrary initial conditions, then clearly $\mathcal{R}%
=0$ on $M\times M$ by (15). Conversely, if $\mathcal{R}=0$ on $M\times M,$
then according to our theorem, any initial condition (= any $1$-arrow) of $%
\mathcal{N}_{1}(w)$ integrates uniquely to a local solution of (13)
satisfying this initial condition. In this case we call $\mathcal{N}_{1}(w)$
uniquely locally integrable. Using a symbolic notation without the indices,
we rewrite (17) as

\begin{equation}
\mathcal{R(}x,y)=I(w;y)-(x,y,f_{1})_{\ast }I(w;x)
\end{equation}

According to (18), we translate $I(w;x)$ from $(U,x)$ to $(V,y)$ using the $%
1 $-arrow $(x,y,f_{1})$ of $\mathcal{N}_{1}(w)$ and compare the value $%
(x,y,f_{1})_{\ast }I(w;x)$ on $(V,y)$ with $I(w;y)$ by subtracting and $%
\mathcal{R(}x,y)$ measures the difference. Consequently, $\mathcal{R=}0$ if
and only if $I(w)$ is invariant with respect to $\mathcal{N}_{1}(w).$
Therefore we can state

\begin{proposition}
The following are equivalent for the subgroupoid $\mathcal{N}_{1}(w)$ $%
\subset \mathcal{U}_{1}.$

1) $\mathcal{R=}0$ on $M\times M$

2) $\mathcal{N}_{1}(w)$ is uniquely locally integrable,i.e, any $1$-arrow of 
$\mathcal{N}_{1}(w)$ integrates uniquely to a local solution of (13).

3) $\mathcal{N}_{1}(w),$ which leaves $w$ invariant by its definition,
leaves also $I(w)$ invariant
\end{proposition}

It is worthwhile here to take a closer look at $\mathcal{R}$ defined by
(17), (18). Clearly $\mathcal{R}(x,x)=0,$ i.e., $\mathcal{R}$ vanishes on
the diagonal of $\mathcal{U}_{0}=M\times M.$ Now $\mathcal{R}$ is a $2$-form
on $M$ and assigns to any $(\overline{x},\overline{y})\in \mathcal{U}_{0}$ a
vector in the fiber of $T\rightarrow M$ over $\overline{y},$ i.e., a tangent
vector at $\overline{y}.$ In the language of bisections, the $2$-form $%
\mathcal{R}$ maps the bisection $(x,f(x),f_{1}(x))$ of $\mathcal{N}_{1}(w)$
to a section of the vector bundle $T\rightarrow M$ over $f(x).$ Therefore,

\begin{equation}
\mathcal{R}:\mathcal{U}_{0}=M\times M\longrightarrow \wedge ^{2}(T^{\ast
})\otimes T
\end{equation}

In (19), we use the same notation for the bundles and the (bi)sections of
these bundles. Note that $\mathcal{R}$ lifts the bisection $(x,f(x))\in 
\mathcal{U}_{0}$ (same abuse of notation) to $(x,f(x),f_{1}(x))\in $ $%
\mathcal{N}_{1}(w)$ by $\varepsilon $ and then maps it to a section of $%
\wedge ^{2}(T^{\ast })\otimes T$ over $f(x).$

Now our purpose is to linearize the above nonlinear picture. Consider the
vector bundle $J_{1}T\rightarrow M$ of $1$-jets of vector fields and let $%
\xi =(\xi ^{i}(x))$ be a vector field integrating to the $1$-parameter
subgroup $f^{i}(t,x)$ of local diffeomorphisms. Therefore

\begin{equation}
f^{i}(0,x)=x^{i}\text{ \ \ \ \ \ \ \ \ \ \ \ }\left[ \frac{\partial
f^{i}(t,x)}{\partial t}\right] _{t=0}=\xi ^{i}(x)
\end{equation}

and

\begin{equation}
\frac{\partial }{\partial x^{j}}\left[ \frac{\partial f^{i}(t,x)}{\partial t}%
\right] _{t=0}=\frac{\partial \xi ^{i}(x)}{\partial x^{j}}=\frac{\partial }{%
\partial t}\left[ \frac{\partial f^{i}(t,x)}{\partial x^{j}}\right] _{t=0}
\end{equation}

\begin{definition}
A (smooth) path in $\mathcal{U}_{1}$ through the identity at $\overline{x}%
\in M$ consists of

1) A path $x(t)$ in $M,$ $t\in (-\epsilon ,\epsilon ),$ $x(0)=\overline{x}$

2) A (smooth) choice of $1$-arrows of $\mathcal{U}_{1}$ from $\overline{x}%
=x(0)$ to $x(t)$ which is identity for $t=0$
\end{definition}

In coordinates, such a path is of the form

\begin{equation}
(\overline{x}^{i},x^{i}(t),f_{j}^{i}(t))\text{ \ \ \ \ }t\in (-\epsilon
,\epsilon )
\end{equation}

\begin{equation*}
(\overline{x}^{i},x^{i}(0),f_{j}^{i}(0))=(\overline{x}^{i},\overline{x}%
^{i},\delta _{j}^{i})
\end{equation*}

Given such a path in $\mathcal{U}_{1}$ (which we will simply call a path),
we define its tangent at $t=0$ by

\begin{equation}
\frac{d}{dt}\left[ (\overline{x}^{i},x^{i}(t),f_{j}^{i}(t))\right] _{t=0}%
\overset{def}{=}(\overline{x}^{i},\xi ^{i}(\overline{x}),\xi _{j}^{i}(%
\overline{x}))
\end{equation}%
and (20), (21) show that $(\overline{x}^{i},\xi ^{i}(\overline{x}),\xi
_{j}^{i}(\overline{x}))$ (briefly $(\xi ^{i}(\overline{x}),\xi _{j}^{i}(%
\overline{x}))$ is an element in the fiber of $J_{1}T\rightarrow M$ over $%
\overline{x}.$ Using (21), we see that all vectors in the fiber over $%
\overline{x}$ are obtained in this way. Indeed, given $(\xi ^{i}(\overline{x}%
),\xi _{j}^{i}(\overline{x})),$ we choose \textit{any} vector field $\xi
=(\xi ^{i}(x))$ around $\overline{x}$ satisfying $\frac{\partial \xi ^{i}(%
\overline{x})}{\partial x^{j}}=\xi _{j}^{i}(\overline{x}),$ and $(\overline{x%
}^{i},f^{i}(t,\overline{x}),\frac{\partial f^{i}(t,\overline{x})}{\partial
x^{j}})$ has the tangent $(\xi ^{i}(\overline{x}),\xi _{j}^{i}(\overline{x}%
)) $ by (20), (21). Henceforth we will use the notations

\begin{eqnarray}
\xi _{0}(\overline{x}) &=&(\xi ^{i}(\overline{x}))\text{ = a tangent vector
at }\overline{x} \\
\xi _{1}(\overline{x}) &=&(\xi ^{i}(\overline{x}),\xi _{j}^{i}(\overline{x}))%
\text{ = a }1\text{-jet of a vector field at }\overline{x}  \notag
\end{eqnarray}

So $\xi _{0}(\overline{x})$ is a vector in the fiber of $J_{0}T=T\rightarrow
M$ over $\overline{x},$ and $\xi _{1}(\overline{x})$ is a vector in the
fiber of $J_{1}T\rightarrow M$ over $\overline{x},$ that projects to $\xi
_{0}(\overline{x}).$ We call this process of passing from the groupoid $%
\mathcal{U}_{1}$ to the vector bundle $J_{1}T\rightarrow M$ (which is
actually the algebroid of $\mathcal{U}_{1})$ linearization. Now our purpose
is to linearize the subgroupoid $\mathcal{N}_{1}(w)\subset \mathcal{U}_{1}$
in the same way whose linearization $N_{1}(w)\rightarrow M$ will be a
subbundle of $J_{1}T\rightarrow M.$

A path in $\mathcal{N}_{1}(w)\subset \mathcal{U}_{1}$ is of the form $(%
\overline{x},x(t),\varepsilon (\overline{x},x(t))$ with the tangent

\begin{eqnarray}
(\xi ^{i}(\overline{x}),\xi _{j}^{i}(\overline{x})) &=&\frac{d}{dt}\left[
x^{i}(t),\varepsilon _{j}^{i}(\overline{x},x(t))\right] _{t=0}=(\xi ^{i}(%
\overline{x}),\Gamma _{aj}^{i}(\overline{x})\xi ^{a}(\overline{x}))  \notag
\\
&=&(\xi ^{i}(\overline{x}),\frac{\partial w_{b}^{i}(\overline{x})}{\partial
x^{a}}\widetilde{w}_{j}^{b}(\overline{x})\xi ^{a}(\overline{x}))
\end{eqnarray}

Using (25) we define the fiber $N_{1}(w)^{\overline{x}}$ $\subset \left(
J_{1}T\right) ^{\overline{x}}$ as those tangents $(\xi ^{i}(\overline{x}%
),\xi _{j}^{i}(\overline{x}))$ satisfying

\begin{equation}
\xi _{j}^{i}(\overline{x})=\Gamma _{aj}^{i}(\overline{x})\xi ^{a}(\overline{x%
})
\end{equation}%
which is clearly a subspace. We define the bundle of vectors $N_{1}(w)%
\overset{def}{=}\cup _{\overline{x}\in M}N_{1}(w)^{\overline{x}}$ with the
obvious projection $N_{1}(w)\rightarrow M$ which is easily seen to be a
vector subbundle. Note that (26) gives a splitting $\varepsilon $ (using the
same notation) of the projection

\begin{equation}
0\longrightarrow T\otimes T^{\ast }\longrightarrow J_{1}(T)\longrightarrow
T\longrightarrow 0
\end{equation}%
defined by

\begin{equation}
\varepsilon :(\xi ^{i})\longrightarrow (\xi ^{i},\Gamma _{aj}^{i}\xi ^{a})
\end{equation}

The process of linearization is defined for all PHG's (more generally for
all Lie groupoids) yielding their algebroids. Since the linear connection
(28) (denoted by $\nabla $ in [Or1]) drops out of this process, it will not
be of primary importance for us. Put more succintly, the splitting (28) is
the linearization of the above nonlinear splitting $\varepsilon $ which is 
\textit{not a connection }for otherwise the curvature of $\varepsilon $
would always vanish since $\varepsilon $ is a trivialization of the
principle bundle $\mathcal{U}_{1}^{e,\bullet }\rightarrow M$ whereas just
the opposite is true for $\mathcal{R}:$ It is surely not always zero and a
very subtle object!

To summarize, a section of $N_{1}(w)\rightarrow M$ is of the form $(\xi
^{i}(x),\Gamma _{aj}^{i}(x)\xi ^{a}(x)).$

\begin{definition}
The section $(\xi ^{i}(x),\xi _{j}^{i}(x))$ of $J_{1}(T)\rightarrow M$ is
prolonged (or holonomic) if

\begin{equation}
\xi _{j}^{i}(x)=\frac{\partial \xi ^{i}(x)}{\partial x^{j}}
\end{equation}
\end{definition}

So a vector field $\xi _{0}=\xi =(\xi ^{i})$ defines a section of $%
J_{1}(T)\rightarrow M$ by prolonging as $pr(\xi )\overset{def}{=}(\xi ^{i},%
\frac{\partial \xi ^{i}}{\partial x^{j}}).$ Therefore, the section $(\xi
^{i}(x),\Gamma _{aj}^{i}(x)\xi ^{a}(x))$ of $N_{1}(w)\rightarrow M$ is
prolonged if and only if

\begin{equation}
\frac{\partial \xi ^{i}(x)}{\partial x^{j}}=\Gamma _{aj}^{i}(x)\xi
^{a}(x)=w_{b}^{i}(y)\widetilde{w}_{a}^{b}(x)\xi ^{a}(x)
\end{equation}

Now (30) is a first order linear system and has unique local solutions with
arbitrary initial conditions if and only if its integrability conditions are
identically satisfied. Some computation shows that these integrability
conditions are given by

\begin{equation}
\mathfrak{R}_{kj,a}^{i}\xi ^{a}\overset{def}{=}\left[ \frac{\partial \Gamma
_{aj}^{i}}{\partial x^{k}}+\Gamma _{bj}^{i}\Gamma _{ak}^{b}\right]
_{[kj]}\xi ^{a}=0
\end{equation}%
for all $\xi =(\xi ^{i}).$ Substituting (7) into (31) expresses $\mathfrak{R}
$ in terms of the structure object $w.$ Now if $\mathfrak{R}=0$ near $%
\overline{x},$ then for any tangent in the fiber $N_{1}(w)^{\overline{x}}$
as initial condition, there is a unique vector field $\xi $ near $\overline{x%
}$ which solves (31) and satisfies the given initial condition. If $%
\mathfrak{R}=0$ on $M,$ we call $N_{1}(w)$ uniquely locally integrable.
Therefore

\begin{equation}
\mathfrak{R=}0\iff N_{1}(w)\text{ is uniquely locally integrable}
\end{equation}

\begin{definition}
$\mathfrak{R}$ is the linear curvature of the subgroupoid $\mathcal{N}%
_{1}(w)\subset \mathcal{U}_{1}.$
\end{definition}

Note that $\mathfrak{R}(x)$ is a $2$-form at $x$ which assigns to the
tangent vector $(\xi ^{i}(x))$ the tangent vector $\mathfrak{R}_{\bullet
,a}^{i}(x)\xi (x)^{a}$ at $x.$ Equivalently, we have

\begin{equation}
\mathfrak{R}:T\longrightarrow \wedge ^{2}(T^{\ast })\otimes T
\end{equation}

Yet another interpretation is that $\mathfrak{R}$ is a $2$-form on $M$ with
values in the vector bundle $Hom(T,T)\rightarrow M.$ At this point,
inspecting (19) and (33) carefully and noting that $T\rightarrow M$ is the
linearization (the algebroid) of $\mathcal{U}_{0}=M\times M$ (the pair
groupoid), it is natural to expect that $\mathfrak{R}$ will be in some sense
the "linearization" of $\mathcal{R}.$

(32) is the linear analog of the equivalence $1)\iff 2)$ in Proposition 7
and it remains to find the linearization of $3)$ of Proposition 7. So let $(%
\overline{x}^{i},x^{i}(t),f_{j}^{i}(t))$ be a path in $\mathcal{U}_{1}$ with
the tangent $\xi _{1}(\overline{x})$ and $\alpha (t)$ a first order
geometric object, for instance a tensor field, defined on this path. Our
purpose is to define the change of $\alpha $ at $\overline{x}$ in the
direction of $\xi _{1}(\overline{x})$ which we will denote by $\left( 
\mathcal{L}_{\{\xi _{1}(\overline{x})\}}\alpha \right) (\overline{x}).$ Note
that $\left( \mathcal{L}_{\{\xi _{1}(\overline{x})\}}\alpha \right) (%
\overline{x})$ should not depend on the path but only on its tangent at $%
\overline{x}.$ When doing this, we should keep in mind the definition of the
ordinary Lie derivative $\mathcal{L}_{\xi _{0}}\alpha $ of $\alpha $ with
respect to some vector field $\xi _{0}$ where both $\xi _{0}$ and $\alpha $
are defined in some neighborhood of $\overline{x}.$ Now the path $(\overline{%
x}^{i},x^{i}(t),f_{j}^{i}(t))$ maps the value $\alpha (\overline{x})=\alpha
(x(0))$ to $(\overline{x}^{i},x^{i}(t),f_{j}^{i}(t))_{\ast }(\alpha (%
\overline{x}))$ and the idea is to compare $(\overline{x}%
^{i},x^{i}(t),f_{j}^{i}(t))_{\ast }(\alpha (\overline{x}))$ with $\alpha
(x(t))$ by dividing their difference by $t$ and letting $t\rightarrow 0.$
The quantity $(\overline{x}^{i},x^{i}(t),f_{j}^{i}(t))_{\ast }(\alpha (%
\overline{x}))$ is computed using the transformation rule of the tensor $%
\alpha .$ For instance, let $\alpha =(\alpha _{j}^{i})$ be a (1,1)-tensor
field which transforms according to

\begin{equation}
\alpha _{j}^{i}(y)=\alpha _{b}^{a}(x)\frac{\partial y^{i}}{\partial x^{a}}%
\frac{\partial x^{b}}{\partial y^{j}}
\end{equation}

Now (34) states the invariance of $\alpha $ with respect to $f,$ i.e., the
local diffeomorphism $f:(U,x)\rightarrow (V,y),$ $y=f(x),$ maps $\alpha (x)$
to $f_{\ast }\alpha (x)$ defined by the RHS of (34) which need not be equal
to $\alpha (y)$ on the LHS of (34) unless $\alpha $ is left invariant by $f.$
Now let $\xi =(\xi ^{i}(x))$ be any vector field defined near $\overline{x}$
satisfying $j_{1}(\xi )(\overline{x})=\xi _{1}(\overline{x}),$ i.e., $\frac{%
\partial \xi ^{i}(\overline{x})}{\partial x^{j}}=\xi _{j}^{i}(\overline{x})$
and let $f^{i}(t,x)$ be the 1-parameter local diffeomorphisms defined by $%
\xi $ with $f^{-1}(t,x)=g(t,y).$ As we observed above, $(\overline{x}%
^{i},f^{i}(t,\overline{x}),\frac{\partial f^{i}(t,\overline{x})}{\partial
x^{j}})$ is a path (possibly different from the above one) having the same
tangent $\xi _{1}(\overline{x})=(\xi ^{i}(\overline{x}),\xi _{j}^{i}(%
\overline{x}))$ at $\overline{x}.$ We define

\begin{equation*}
\left( \mathcal{L}_{\{\xi _{1}(\overline{x})\}}\alpha \right) _{j}^{i}(%
\overline{x})\overset{def}{=}
\end{equation*}

\begin{eqnarray}
&&\frac{d}{dt}\left[ \alpha _{j}^{i}(f(t,\overline{x}))-\alpha _{b}^{a}(%
\overline{x})\frac{\partial f^{i}(t,\overline{x})}{\partial x^{a}}\frac{%
\partial g^{b}(t,\overline{y})}{\partial y^{j}}\right] _{t=0} \\
&=&\left[ \frac{\partial \alpha _{j}^{i}(f(t,\overline{x}))}{\partial y^{a}}%
\frac{df^{a}(t,\overline{x})}{dt}\right] _{t=0}-\alpha _{b}^{a}(\overline{x})%
\frac{d}{dt}\left[ \frac{\partial f^{i}(t,\overline{x})}{\partial x^{a}}%
\right] _{t=0}\left[ \frac{\partial g^{b}(t,\overline{y})}{\partial y^{j}}%
\right] _{t=0}  \notag \\
&&-\alpha _{b}^{a}(\overline{x})\left[ \frac{\partial f^{i}(t,\overline{x})}{%
\partial x^{a}}\right] _{t=0}\frac{d}{dt}\left[ \frac{\partial g^{b}(t,%
\overline{y})}{\partial y^{j}}\right] _{t=0}  \notag \\
&=&\frac{\partial \alpha _{j}^{i}(\overline{x})}{\partial x^{a}}\xi ^{a}(%
\overline{x})-\alpha _{b}^{a}(\overline{x})\frac{\partial \xi ^{i}(\overline{%
x})}{\partial x^{a}}(\delta _{j}^{b})-\alpha _{b}^{a}(\overline{x})(\delta
_{a}^{i})(-\frac{\partial \xi ^{b}(\overline{x})}{\partial x^{j}}) \\
&=&\frac{\partial \alpha _{j}^{i}(\overline{x})}{\partial x^{a}}\xi ^{a}(%
\overline{x})-\alpha _{j}^{a}(\overline{x})\xi _{a}^{i}(\overline{x})+\alpha
_{a}^{i}(\overline{x})\xi _{j}^{a}(\overline{x})
\end{eqnarray}

We make six important observations.

1) To compute $\left( \mathcal{L}_{\{\xi _{1}(\overline{x})\}}\alpha \right)
(\overline{x}),$ (35) shows that we need only the values of $\alpha $ on the
path and (37) shows that $\left( \mathcal{L}_{\{\xi _{1}(\overline{x}%
)\}}\alpha \right) (\overline{x})$ depends only on the tangent of the path
at $\overline{x}$ as required. Furthermore, $\mathcal{L}_{\{\xi _{1}(%
\overline{x})\}}\alpha $ is linear in the argument $\xi _{1}.$

2) We gave above all the details in the derivation of (37). For
computational purposes however, all we have to do to derive (37) is to
formally substitute $y^{i}=x^{i}+t\xi ^{i}(x),$ $\frac{\partial y^{i}}{%
\partial x^{j}}=\delta _{j}^{i}+t\frac{\partial \xi ^{i}}{\partial x^{j}},$ $%
x^{i}=y^{i}-t\xi ^{i}(y),$ $\frac{\partial x^{i}}{\partial y^{j}}=\delta
_{j}^{i}-t\frac{\partial \xi ^{i}}{\partial y^{j}}$ into (34), collect all
the terms on the LHS of (34) and differentiate the resulting expression with
respect to $t$ at $t=0.$

3) $\left( \mathcal{L}_{\{pr(\xi )\}}\alpha \right) (x)$ is the ordinary Lie
derivative $\mathcal{L}_{\xi }\alpha $ of $\alpha $ with respect to $\xi $
as expected. Therefore, (36) is the ordinary Lie derivative of $\alpha $ and
the passage from (36) to (37) shows that $\mathcal{L}_{\{\xi _{1}\}}\alpha $
is computed by first computing the ordinary Lie derivative of $\alpha $ and
replacing the derivatives of the vector field $\xi $ with jet variables.

4) If $\alpha $ is a global (1,1)-tensor field on $M$ and $\xi _{1}$ is a
global section of $J_{1}T\rightarrow M,$ then $\mathcal{L}_{\{\xi
_{1}\}}\alpha $ is defined pointwise by (37) and is another tensor of the
same type as $\alpha .$

5) $\mathcal{L}_{\{\xi _{1}\}}$ is a derivation of the tensor algebra and
commutes with contractions.

6) $\mathcal{L}_{\{\xi _{1}\}}\alpha $ is defined in the same way for all
first order geometric objects $\alpha ,$ in particular, for all tensor
fields.

\begin{definition}
$\mathcal{L}_{\{\xi _{1}\}}\alpha $ is the formal Lie derivative of $\alpha $
with respect to the section $\xi _{1}$ of $J_{1}T\rightarrow M.$
\end{definition}

Now let us choose our path in the definition of $\mathcal{L}_{\{\xi
_{1}\}}\alpha $ in $\mathcal{N}_{1}(w)\subset \mathcal{U}_{1}$ keeping in
mind that $\xi _{1}=\varepsilon \xi _{0}=\varepsilon \xi $ in this case.
Substituting (28) into (37), we find

\begin{eqnarray}
\left( \mathcal{L}_{\varepsilon \xi }\alpha \right) _{j}^{i}(\overline{x})
&=&\frac{\partial \alpha _{j}^{i}(\overline{x})}{\partial x^{a}}\xi ^{a}(%
\overline{x})-\alpha _{j}^{b}(\overline{x})\Gamma _{ab}^{i}(\overline{x})\xi
^{a}(\overline{x})+\alpha _{b}^{i}(\overline{x})\Gamma _{aj}^{b}(\overline{x}%
)\xi ^{a}(\overline{x})  \notag \\
&=&\left( \frac{\partial \alpha _{j}^{i}(\overline{x})}{\partial x^{a}}%
-\Gamma _{ab}^{i}(\overline{x})\alpha _{j}^{b}(\overline{x})+\Gamma
_{aj}^{b}(\overline{x})\alpha _{b}^{i}(\overline{x})\right) \xi ^{a}(%
\overline{x})
\end{eqnarray}

We observe that another linear connection pops up from (38) (denoted by $%
\widetilde{\nabla }$ in [Or1]) and the remarkable fact is that this
connection \textit{differs from the above one by the integrability object }$%
I(w)_{jk}^{i}$ \textit{called torsion in [Or1]!} \textit{Furthermore, this
new connection has vanishing curvature! }On the hand, the formal Lie
derivative $\mathcal{L}_{\xi _{k+1}}\alpha _{k}$ is defined in the obvious
way for all $k^{\prime }$th order geometric object $\alpha _{k}$ and for a
section $\xi _{k+1}$ of $J_{k+1}(T)\rightarrow M$ by computing the ordinary
Lie derivative $\mathcal{L}_{\xi }\alpha _{k}$ and replacing the derivatives
of $\xi $ with jet variables. For a PHG of order $k,$ $\mathcal{L}_{\xi
_{k}}\alpha _{k}$ is defined by $\mathcal{L}_{\varepsilon \xi _{k}}\alpha
_{k}.$ We observe the crucial fact that $\mathcal{L}_{\varepsilon \xi
_{k}}\alpha _{k}$ is \textit{not a connection any more for }$k\geq 2$ but
forms the basis of an abstraction called "algebroid connection" in the
modern theory of algebroids.

It is easy to check that $\mathcal{L}_{\varepsilon \xi }w=0$ for all $\xi .$
Indeed, $\mathcal{L}_{\varepsilon \xi }\alpha $ measures how $\alpha $
changes along the paths of $\mathcal{N}_{1}(w)$ and $w$ is constant along
all these paths by the definition of $\mathcal{N}_{1}(w).$ Now we claim that 
$\mathcal{L}_{\varepsilon \xi }I(w)=\mathfrak{R(\xi )},$ or in more detail

\begin{equation}
\left[ \frac{d}{dt}\mathcal{R}_{kj}^{i}\mathcal{(}x,x+t\xi )\right] _{t=0}=%
\mathcal{L}_{\varepsilon \xi }\left( I(w;x)_{kj}^{i}\right) =\mathfrak{R}%
_{kj,a}^{i}(x)\xi ^{a}
\end{equation}

The first equality in (39) holds by the definition of $\mathcal{R}$ by (17)
and it justifies, in view of the second equality in (39), that $\mathfrak{R}$
is indeed the linearization of $\mathcal{R}$ as forseen above. To prove the
second equality, all we need to observe is that $\mathcal{L}_{\varepsilon
\xi }=\widetilde{\nabla }_{\xi }$ and $I(w)_{kj}^{i}=T_{kj}^{i}$ in [Or1]
and (39) is in fact the \textit{definition of }$\mathfrak{R}$ \textit{in
[Or1] }(see Definition 6.1)\textit{. }Therefore, the linearization of $3)$
of Proposition 7 that we search for is given by the middle term in (39)!
Combining (39) with (32), we now state

\begin{proposition}
The following are equivalent.

1) $\mathfrak{R}=0$

2) $N_{1}(w)\rightarrow M$ is uniquely locally integrable

3) $\mathcal{L}_{\varepsilon \xi }I(w)=0$ for all $\xi \in T$
\end{proposition}

At this stage, it is easy to guess that $3)$ of Proposition 7 and $3)$ of
Proposition 12 are equivalent. Indeed, in the language of [Or1], $3)$ of
Proposition 7 asserts that $I(w)$ is $\varepsilon $-invariant and $3)$ of
Proposition 12 asserts $I(w)$ is $\widetilde{\nabla }$-parallel and these
two concepts are equivalent according to Proposition 5.5 in [Or1]. A more
intuitive and amusing argument goes as follows: Clearly $\mathcal{R}=0$
implies $\mathfrak{R}=0$ by (39). Conversely, if $\mathfrak{R}=0,$ then the
"derivative of $\mathcal{R}(x,y)$ with respect $y"$ vanishes identically
according to (39) and therefore $\mathcal{R}(x,y)$ is "constant in $y".$
Therefore $\mathcal{R}=0$ on $M\times M$ since $\mathcal{R}(x,x)=0.$

Thus we state

\begin{proposition}
(Lie's 3'rd Theorem) The conditions of Proposition 5 and Proposition 12 are
equivalent.
\end{proposition}

It is explained in [Or1] why Proposition 13 is called Lie's 3'rd Theorem.

There is another fundamental operator lurking in the above picture and this
is a good place to pinpoint it. Consider the Spencer operator

\begin{eqnarray}
D &:&J_{1}(T)\longrightarrow \wedge (T^{\ast })\otimes T  \notag \\
&:&(\xi ^{i},\xi _{j}^{i})\longrightarrow (\frac{\partial \xi ^{i}}{\partial
x^{j}}-\xi _{j}^{i})
\end{eqnarray}%
which restricts to

\begin{eqnarray}
D &:&N_{1}(w)\longrightarrow \wedge (T^{\ast })\otimes T \\
&:&(\xi ^{i},\Gamma _{aj}^{i}\xi ^{a})\longrightarrow (\frac{\partial \xi
^{i}}{\partial x^{j}}-\Gamma _{aj}^{i}\xi ^{a})  \notag
\end{eqnarray}

Therefore $D_{\eta }\xi \overset{def}{=}\left( \frac{\partial \xi ^{i}}{%
\partial x^{b}}-\Gamma _{ab}^{i}\xi ^{a}\right) \eta ^{b}=\nabla _{\eta }\xi
,$ i.e., $D_{\eta }=\nabla _{\eta }$ (see (28)). As we remarked above, $%
\mathcal{L}_{\{\eta _{k}\}}$ will not be a connection on the PHG for $k\geq
1 $ and we will see below that $D_{\eta _{k}}$ will sometimes become a "very
special connection" but under a strong assumption!

\section{Affine PHG's}

In this section we will imitatate our above arguments \textit{word by word}
and therefore will not give all the details but eloborate when a new
phenomenon occurs.

Consider the transformation group $Aff(\mathbb{R}^{n})=GL(n,\mathbb{R}%
)\ltimes \mathbb{R}^{n}$ of $\mathbb{R}^{n}.$ Fixing, for instance, the
origin $o\in \mathbb{R}^{n}$ and its stabilizer $GL(n,\mathbb{R)},$ we
define the map $j_{k}:GL(n,\mathbb{R)}=G_{1}(n)\mathbb{\rightarrow }G_{k}(n)$
by $g\rightarrow j_{k}(g)^{o}$ which is injective for $k\geq 1,$ where $%
G_{k}(n)$ is the $k^{\prime }$th order jet group in $n$ variables.
Identifying $G_{1}(n)$ with its image $j_{2}(G_{1}(n))\subset G_{2}(n),$ the
left coset space $G_{2}(n)/G_{1}(n))$ is parametrized by functions $(\Gamma
_{jk}^{i})$ (see [Or1], 175-178 for details). Therefore, the left action of $%
G_{2}(n)$ on $G_{2}(n)/G_{1}(n))$ defines a geometric object $\Gamma $ with
components $(\Gamma _{jk}^{i}(x))$ subject to the transformation rule

\begin{equation}
\Gamma _{ab}^{i}(f(x))\frac{\partial f^{a}(x)}{\partial x^{j}}\frac{\partial
f(x)^{b}}{\partial x^{k}}=\Gamma _{jk}^{a}(x)\frac{\partial f^{i}(x)}{%
\partial x^{a}}+\frac{\partial ^{2}f^{a}(x)}{\partial x^{j}\partial x^{k}}
\end{equation}%
upon a coordinate change $(x)\rightarrow (f(x)).$ It is standard to
interpret $\Gamma =(\Gamma _{jk}^{i}(x))$ as a "torsionfree affine
connection" but we will carefully avoid this interpretation and regard $%
\Gamma $ merely as a geometric object on $M$ like $w.$Using (42), we
consider those $2$-arrows $(\overline{x},\overline{y},\overline{f}_{1},%
\overline{f}_{2})=(\overline{x}^{i},\overline{y}^{i},\overline{f}_{j}^{i},%
\overline{f}_{jk}^{i})$ of $\mathcal{U}_{2}$ which preserve $\Gamma ,$ that
is

\begin{equation}
\Gamma _{ab}^{i}(\overline{y})\overline{f}_{j}^{a}\overline{f}%
_{k}^{b}=\Gamma _{jk}^{a}(\overline{x})\overline{f}_{a}^{i}+\overline{f}%
_{jk}^{i}
\end{equation}

\begin{definition}
The subgroupoid $\mathcal{H}_{2}(\Gamma )\subset \mathcal{U}_{2}$ defined by
(43) is an affine PHG on $M$ and $\Gamma $ is its structure object.
\end{definition}

We observe that $(\overline{x}^{i},\overline{y}^{i},\overline{f}_{j}^{i})\in 
\mathcal{U}_{1}$ is arbitrary in (43) and $\overline{f}_{jk}^{i}$ is
uniquely determined by $(\overline{x}^{i},\overline{y}^{i},\overline{f}%
_{j}^{i}),$ which gives the splitting

\begin{equation*}
\varepsilon :\mathcal{H}_{1}(\Gamma )=\mathcal{U}_{1}\longrightarrow 
\mathcal{U}_{2}
\end{equation*}%
\begin{equation}
\varepsilon :(\overline{x}^{i},\overline{y}^{i},\overline{f}%
_{j}^{i})\longrightarrow (\overline{x}^{i},\overline{y}^{i},\overline{f}%
_{j}^{i},\text{ }\Gamma _{ab}^{i}(\overline{y})\overline{f}_{j}^{a}\overline{%
f}_{k}^{b}-\Gamma _{jk}^{a}(\overline{x})\overline{f}_{a}^{i})
\end{equation}

If we set $\overline{x}=\overline{y}$ in (43), we get the the vertex groups $%
\mathcal{H}_{2}(\Gamma )^{\overline{x}}\cong \mathcal{H}_{1}(\Gamma )^{%
\overline{x}}$ at $\overline{x}.$ We can always find coordinates $(x)$
around $\overline{x}$ with the property $\Gamma _{jk}^{i}(\overline{x})=0$
as can be seen from (43). We call $(x)$ regular coordinates at $\overline{x}%
. $ In such coordinates, the vertex groups are identified with $%
G_{1}(n)=GL(n,\mathbb{R}).$

We can clearly replace the arrows in (43) and (44) by bisections as before.
A bisection $(x^{i},f^{i}(x),f_{j}^{i}(x),f_{jk}^{i}(x))$ of $\mathcal{U}%
_{2} $ is prolonged if $\frac{\partial f^{i}(x)}{\partial x^{j}}%
=f_{j}^{i}(x) $ and $\frac{\partial f_{j}^{i}(x)}{\partial x^{k}}%
=f_{jk}^{i}(x).$ Therefore, a bisection of $\mathcal{H}_{2}(\Gamma )$ is
prolonged if and only if

\begin{eqnarray}
\frac{\partial f^{i}(x)}{\partial x^{j}} &=&f_{j}^{i}(x) \\
\frac{\partial f_{j}^{i}(x)}{\partial x^{k}} &=&f_{jk}^{i}(x)=\Gamma
_{ab}^{i}(f(x))f_{j}^{a}(x)f_{k}^{b}(x)-\Gamma _{jk}^{a}(x)f_{a}^{i}(x) 
\notag
\end{eqnarray}

Clearly (45) has a solution if and only if (42) has a solution. However,
(42) is a second order PDE whereas (45) is a \textit{first order system }of
PDE's. Note that (45) is a \textit{closed }system in the sense that it
expresses the derivatives of the unknown functions $f^{i}(x),f_{j}^{i}(x)$
in terms of themselves and $x$'s due to the splitting $\varepsilon .$
Therefore, we see that a prolonged bisection corresponds to the well known
trick of introducing jet variables to reduce a second order PDE to a first
order system. Now (44) serves as the initial conditions for (45). The
integrability conditions of (45) are given by

\begin{eqnarray}
\frac{\partial f_{j}^{i}}{\partial x^{k}}-\frac{\partial f_{k}^{i}}{\partial
x^{j}} &=&0 \\
\frac{\partial f_{jk}^{i}}{\partial x^{r}}-\frac{\partial f_{rk}^{i}}{%
\partial x^{j}} &=&0  \notag
\end{eqnarray}

Since the second expression of (45) is symmetric in $j,k,$ the first
condition of (46) is identically satisfied. To check the second, we
differentiate the second expression of (45) with respect to $x^{r},$
substitute back from (45) and alternate $r,j.$ After some straightforward
computation, we find

\begin{eqnarray}
\mathcal{R}_{rj,k}^{i}(x,f(x),f_{1}(x))\overset{def}{=} &&\left[ \frac{%
\partial \Gamma _{jk}^{i}(f(x))}{\partial y^{r}}+\Gamma
_{jb}^{i}(f(x))g_{rk}^{b}(f(x))\right] _{[rj]} \\
&&-g_{j}^{a}g_{r}^{c}g_{k}^{d}\left[ \frac{\partial \Gamma _{ad}^{e}(x)}{%
\partial x^{c}}+\Gamma _{ab}^{e}(x)\Gamma _{cd}^{b}(x)\right]
_{[ca]}f_{e}^{i}(x)=0  \notag
\end{eqnarray}%
where $(x,f(x),f_{1}(x))=(x^{i},f^{i}(x),f_{j}^{i}(x))$ is a bisection of $%
\mathcal{H}_{1}(\Gamma )=\mathcal{U}_{1},$ $g_{1}=f_{1}^{-1}.$ Note that $%
\mathcal{R}_{rj,k}^{i}$ is a $2$-form in the indices $r,j$ and a $(1,1)$%
-tensor in the indices\textit{\ }$i,k$

\begin{definition}
$\mathcal{R}$ is the nonlinear curvature of $\mathcal{H}_{2}(\Gamma ).$
\end{definition}

We define 
\begin{equation}
I_{rj,k}^{i}(\Gamma ;x)\overset{def}{=}\left[ \frac{\partial \Gamma
_{jk}^{i}(x)}{\partial x^{r}}+\Gamma _{ja}^{i}(x)\Gamma _{rk}^{a}(x)\right]
_{[rj]}
\end{equation}%
and call $I(\Gamma )=(I_{rj,k}^{i}(\Gamma ;x))$ the integrability object of $%
\Gamma =(\Gamma _{jk}^{i}(x))$ (which is the curvature of the affine
connection $\Gamma !).$ Like (18), we write (47) symbolically as%
\begin{equation}
\mathcal{R(}x,y,f_{1})=I(\Gamma ;y)-(x,y,f_{1})_{\ast }I(\Gamma ;x))=0
\end{equation}%
and (49) expresses the invariance of $I(\Gamma )$ by the bisections (or
equivalently $1$-arrows) of $\mathcal{H}_{1}(\Gamma ).$ Note that $%
f_{1}=(f_{j}^{i})$ is a dependent variable in (18) whereas an independent
variable in (49). We view $\mathcal{R}$ as a map

\begin{equation}
\mathcal{R}:\mathcal{H}_{1}(\Gamma )=\mathcal{U}_{1}\longrightarrow \wedge
^{2}(T^{\ast })\otimes T^{\ast }\otimes T
\end{equation}

We will see in the next section that $\mathcal{R}$ is actually a map

\begin{equation}
\mathcal{R}:\mathcal{H}_{1}(\Gamma )=\mathcal{U}_{1}\longrightarrow \wedge
^{2}(T^{\ast })\otimes J_{1}(T)
\end{equation}%
but its projection on $T$ vanishes (see (27)) since $\mathcal{H}_{1}(\Gamma
)=\mathcal{U}_{1}!$

Therefore, (49) holds for all bisections if and only if all initial
conditions (44) integrate locally and uniquely to local solutions of (42).
Therefore we state

\begin{proposition}
The following are equivalent.

1) $\mathcal{R}=0$ on $\mathcal{U}_{1}$

2) $\mathcal{H}_{2}(\Gamma )$ is uniquely locally integrable

3) $\mathcal{H}_{2}(\Gamma ),$ which leaves $\Gamma $ invariant by
definition, leaves also $I(\Gamma )$ invariant
\end{proposition}

Now a new phenomenon occurs due to do nontriviality of the stabilizers:
Choosing $x=y$ in (49) and noting that $f$ \ is arbitrary, it follows that
the tensor $I(\Gamma (x))$ is fixed by $\mathcal{H}_{2}(\Gamma )^{x}\cong 
\mathcal{H}_{1}(\Gamma )^{x}\cong GL(n,\mathbb{R}).$ This is possible if and
only if $I(\Gamma )=0$ on $M$ and we obtain

\begin{proposition}
The conditions of Proposition 16 are equivalent to $I(\Gamma )=0.$
\end{proposition}

Therefore, an affine PHG is flat if and only if it is flat in the classical
sense, i.e., "the torsionfree affine connection $\Gamma $ has vanishing
curvature". \textit{As a very intriguing fact, however, the linearization of
(50) will be the curvature of a connection on }$J_{1}(T)\rightarrow M$%
\textit{\ and not }$T\rightarrow M$\textit{\ where the affine connection }$%
\Gamma $\textit{\ is defined !!}

The linearization of $\mathcal{H}_{2}(\Gamma )$ is now straightforward. A
path in $\mathcal{U}_{2}$ (through the identity at $\overline{x}\in M)$
consists of a path $x(t)$ in $M,$ $t\in (-\epsilon ,\epsilon ),$ $x(0)=%
\overline{x}$ and a smooth choice of $2$-arrows of $\mathcal{U}_{2}$ from $%
\overline{x}=x(0)$ to $x(t)$ which is identity for $t=0.$ In coordinates,
such a path is of the form $(\overline{x}%
^{i},x^{i}(t),f_{j}^{i}(t),f_{jk}^{i}(t)),$ $t\in (-\epsilon ,\epsilon ),$ $(%
\overline{x}^{i},x^{i}(0),f_{j}^{i}(0),f_{jk}^{i}(0))=(\overline{x}^{i},%
\overline{x}^{i},\delta _{j}^{i},0)$ and the tangent of this path is defined
by

\begin{equation}
\frac{d}{dt}\left[ (\overline{x}^{i},x^{i}(t),f_{j}^{i}(t),f_{jk}^{i}(t))%
\right] _{t=0}=(\xi ^{i}(\overline{x}),\xi _{j}^{i}(\overline{x}),\xi
_{jk}^{i}(\overline{x}))
\end{equation}%
giving all vectors in the fiber of $J_{2}T\rightarrow M$ over $\overline{x}.$
Choosing our paths in $\mathcal{H}_{2}(\Gamma )$ and using (44), we define
the linearization $H_{2}(\Gamma )\rightarrow M$ of the subgroupoid $\mathcal{%
H}_{2}(\Gamma )\subset \mathcal{U}_{2}$ which is a vector subbundle of $%
J_{2}(T)\rightarrow M$ and its fiber over $\overline{x}$ is defined by
(omitting $\overline{x}$ from our notation)

\begin{equation}
\frac{\partial \Gamma _{jk}^{i}}{\partial x^{a}}\xi ^{a}+\Gamma _{ka}^{i}\xi
_{j}^{a}+\Gamma _{ja}^{i}\xi _{k}^{a}-\xi _{a}^{i}\Gamma _{jk}^{a}=\xi
_{jk}^{i}
\end{equation}

Clearly, $H_{0}(\Gamma )=T=$ the tangent bundle of $M,$ $H_{1}(\Gamma
)=J_{1}(T)$ and (53) defines a splitting $\varepsilon :H_{1}(\Gamma
)\rightarrow H_{2}(\Gamma )$ so that $H_{1}(\Gamma )\cong H_{2}(\Gamma ).$
Now (53) defines a second order linear PDE on the vector fields $\xi =(\xi
^{i}(x))$ on $M$ which can be reduced to a first order system by introducing
sections. The integrability conditions of this first order system are easily
obtained from (53) and are of the form

\begin{equation}
\mathfrak{R(\xi }_{1})_{rj,k}^{i}=0
\end{equation}

We will not bother here to give the explicit form of (54). Note that $%
\mathfrak{R(\xi }_{1})$ depends linearly on the section $\xi _{1}$ of $%
J_{1}(T)\rightarrow M.$ Now $\mathfrak{R}$ is a $2$-form on $M$ which maps
sections of $J_{1}(T)\rightarrow M$ linearly to sections of $T^{\ast
}\otimes T\rightarrow M$ (actually to sections of $J_{1}(T)\rightarrow M!),$
i.e.,

\begin{equation}
\mathfrak{R}:J_{1}(T)\longrightarrow \wedge ^{2}(T^{\ast })\otimes T^{\ast
}\otimes T\text{ \ }\subset \text{ \ }\wedge ^{2}(T^{\ast })\otimes J_{1}(T)
\end{equation}%
and can be interpreted also as a $2$-form on $M$ with values in $%
Hom(J_{1},T^{\ast }\otimes T)\subset Hom(J_{1},J_{1}).$ After some
computation (which turns out to be redundant, see below) we deduce $\mathcal{%
L}_{\varepsilon \xi _{1}}I(\Gamma )=\mathfrak{R(\xi }_{1}),$ i.e.,

\begin{equation}
\mathcal{L}_{\varepsilon \xi _{1}}\left( I(\Gamma )_{rj,k}^{i}\right) =%
\mathfrak{R(\xi }_{1})_{rj,k}^{i}
\end{equation}

Therefore we deduce

\begin{proposition}
The following are equivalent.

1) $\mathfrak{R}=0$

2) $H_{2}(\Gamma )\rightarrow M$ is uniquely locally integrable

3) $\mathcal{L}_{\varepsilon \xi _{1}}I(\Gamma )=0$ for all $\xi _{1}\in
H_{1}(\Gamma )=J_{1}(T)$
\end{proposition}

Now sit back a moment and look at what we did so far: We differentiate the
nonlinear equations (45) with respect to $t$ to get the linear system (53)
and we differentiate (53) with respect to $x$ to deduce the linear
integrability conditions (54). However, these two operations commute: We can
differentiate first (45) with respect to $x$ and deduce the nonlinear
integrability conditions (50) and then linearize (50) by differentiating
with respect to $t.$ Hence we conclude that (55) must be the linearization
of (50), that is (omitting the indices) we have

\begin{equation}
\frac{d}{dt}\left[ \mathcal{R}(x,x+t\xi _{(0)},\delta +t\xi _{(1)})\right]
_{t=0}=\mathfrak{R(\xi }_{1})
\end{equation}%
where we used the notation $\xi _{1}=(\xi _{(0)},\xi _{(1)})$ (see (24)) for
the section $\xi _{1}$ of $J_{1}(T)=H_{1}(\Gamma )\rightarrow M.$ Like (39),
it follows from (57) that $\mathcal{R}=0$ implies $\mathfrak{R}=0$ and the
calculus fact suggests that the converse holds too (this can be proved
rigorously by reducing the problem to absolute parallelism, see the next
section).

Thus we state

\begin{proposition}
(Lie's 3'rd Theorem) The conditions of Propositions 16, 17 and 18 are
equivalent.
\end{proposition}

It is worthwhile to emphasize again that our constructions are completely
independent of connections (let alone torsionfree ones), but incorporate
only splittings which are built into the definitions of the geometric
structures and serve the purpose of reducing certain PDE's to first order
systems.

A last remark: As we have seen, what we did for affine PHG's in this section
follows word by word the prescription for absolute parallelism and the same
will be true in the next section for Riemannian PHG's. This fact suggests
the existence of a unique principle that handles all cases at one stroke.
This is indeed true and this process of reduction to parallelism will be
expained in Section 6 below.

\section{Riemannian PHG's}

Replacing $Aff(\mathbb{R}^{n})=GL(n,\mathbb{R})\ltimes \mathbb{R}^{n}$ with $%
Iso(\mathbb{R}^{n})=O(n)\ltimes \mathbb{R}^{n},$ the same argument (see
[Or1] for details) defines a geometric object $\mathbf{g}=(g,\Gamma
)=(g_{ij}(x),\Gamma _{jk}^{i}(x))$ where $g=(g_{ij}(x))$ is a Riemannian
metric and $\Gamma =(\Gamma _{jk}^{i}(x))$ is a "torsionfree affine
connection \textit{not necessarily the Levi-Civita connection}". It is
standard in modern differential geometry to regard the components of $%
\mathbf{g}$ as seperate objects whereas from the present standpoint $\mathbf{%
g}$ is a single object whose components are subject to the transformation
rule

\begin{eqnarray}
g_{ab}(f(x))\frac{\partial f^{a}(x)}{\partial x^{j}}\frac{\partial f^{b}(x)}{%
\partial x^{k}} &=&g_{jk}(x) \\
\Gamma _{ab}^{i}(f(x))\frac{\partial f^{a}(x)}{\partial x^{j}}\frac{\partial
f(x)^{b}}{\partial x^{k}} &=&\Gamma _{jk}^{a}(x)\frac{\partial f^{i}(x)}{%
\partial x^{a}}+\frac{\partial ^{2}f^{a}(x)}{\partial x^{j}\partial x^{k}}
\end{eqnarray}%
upon a coordinate change $(x)\rightarrow (f(x)).$ Clearly (59) is identical
with (42). Using (58)+(59) we consider those $2$-arrows $(\overline{x}^{i},%
\overline{y}^{i},\overline{f}_{j}^{i},\overline{f}_{jk}^{i})$ of $\mathcal{U}%
_{2}$ which preserve $\mathbf{g=}(g,\Gamma ),$ that is

\begin{eqnarray}
g_{ab}(\overline{y})\overline{f}_{j}^{a}\overline{f}_{k}^{b} &=&g_{jk}(%
\overline{x}) \\
\Gamma _{ab}^{i}(\overline{y})\overline{f}_{j}^{a}\overline{f}_{k}^{b}
&=&\Gamma _{jk}^{a}(\overline{x})\overline{f}_{a}^{i}+\overline{f}_{jk}^{i}
\end{eqnarray}

\begin{definition}
The subgroupoid $\mathcal{K}_{2}(\mathbf{g)=}\mathcal{K}_{2}(g,\Gamma
)\subset \mathcal{U}_{2}$ defined by (60)+(61) is a Riemannian PHG on $M.$
\end{definition}

Now (61) defines a splitting

\begin{equation}
\varepsilon :\mathcal{K}_{1}(\mathbf{g)\longrightarrow }\mathcal{K}_{2}(%
\mathbf{g)}
\end{equation}%
which is the restriction of (44) to $\mathcal{K}_{1}(\mathbf{g).}$ We can
always find coordinates $(x)$ around $\overline{x}$ with the property $%
g_{ij}(\overline{x})=\delta _{ij}$ and $\Gamma _{jk}^{i}(\overline{x})=0.$
We call $(x)$ regular coordinates at $\overline{x}.$ In such coordinates,
the vertex group $\mathcal{K}_{2}(\mathbf{g})^{\overline{x}}\cong \mathcal{K}%
_{1}(\mathbf{g})^{\overline{x}}$ is identified with the orthogonal group $%
O(n).$ Clearly we can write (60)+(61) also in terms of bisections.

Substituting $y^{i}=x^{i}+t\xi ^{i},$ $f_{j}^{i}=\delta _{j}^{i}+t\xi
_{j}^{i},$ $f_{jk}^{i}=t\xi _{jk}^{i}$ into (60)+(61) and differentiating
with respect to $t=0,$ we get the defining equations of the linearization $%
K_{2}(\mathbf{g})\rightarrow M$ as a subbundle of $J_{2}(T)\rightarrow M$
given by

\begin{eqnarray}
\frac{\partial g_{jk}}{\partial x^{a}}\xi ^{a}+g_{ka}\xi _{j}^{a}+g_{ja}\xi
_{k}^{a} &=&0 \\
\frac{\partial \Gamma _{jk}^{i}}{\partial x^{a}}\xi ^{a}+\Gamma _{ka}^{i}\xi
_{j}^{a}+\Gamma _{ja}^{i}\xi _{k}^{a}-\xi _{a}^{i}\Gamma _{jk}^{a} &=&\xi
_{jk}^{i}
\end{eqnarray}%
and (64) is clearly identical with (53). Note that for $\xi ^{i}=0,$
(63)+(64) become $\xi _{j}^{k}=-\xi _{k}^{j}$ and $\xi _{jk}^{i}=0$ in
regular coordinates.

Now (63)+(64) define a splitting

\begin{equation*}
\varepsilon :K_{1}(\mathbf{g})\longrightarrow K_{2}(\mathbf{g})
\end{equation*}%
and $\varepsilon $ is the restriction of (44).

As it is clear by now, by introducing bisections, we can reduce the $2$'nd
order nonlinear PDE (60)+(61) to an equivalent first order nonlinear \textit{%
mixed system with constraints }and then compute the integrability conditions 
$\mathcal{R}=0.$ Similarly, by introducing sections, we can reduce the $2$%
'nd order linear PDE (63)+(64) to a first order linear mixed system and
compute the integrability conditions $\mathfrak{R}=0$ keeping in mind that
differentiations with respect to $\ t$ and $x$ commute so that $\mathfrak{R}$
will be the linearization of $\mathcal{R}.$

To pinpoint the new phenomenon, we differentiate

\begin{equation}
g_{ab}(f(x))f_{j}^{a}(x)f_{k}^{b}(x)=g_{jk}(x)
\end{equation}%
with respect to $x^{r}$ and substitute bisections. This gives

\begin{eqnarray}
&&\frac{\partial g_{ab}(f(x))}{\partial y^{c}}%
f_{r}^{c}(x)f_{j}^{a}(x)f_{k}^{b}(x)+g_{ab}(f(x))f_{rj}^{a}(x)f_{k}^{b}(x)+g_{ab}(f(x))f_{j}^{a}(x)f_{rk}^{b}(x)
\notag \\
&=&\frac{\partial g_{jk}(x)}{\partial x^{r}}
\end{eqnarray}

Now we substitute $f_{rk}^{b}(x)$ from (61) into (66) and alternate $r,j$ in
(66). After some straightforward computation, we arrive at

\begin{equation}
\left[ g^{ak}\left( \frac{\partial g_{ja}(y)}{\partial y^{r}}%
+g_{jb}(y)\Gamma _{ra}^{b}(y)\right) \right] _{[rj]}
\end{equation}

\begin{equation*}
-\left[ g^{be}\left( \frac{\partial g_{db}(x)}{\partial x^{c}}%
+g_{da}(x)\Gamma _{cb}^{a}(x)\right) \right]
_{[cd]}f_{r}^{c}(x)f_{j}^{d}(x)f_{e}^{k}(x)=0
\end{equation*}

We define $I_{1}(\mathbf{g};x)$ by

\begin{equation}
I_{1}(\mathbf{g};x)\overset{def}{=}I_{rj}^{k}(\mathbf{g};x)\overset{def}{=}%
\left[ g^{ak}\left( \frac{\partial g_{ja}(x)}{\partial x^{r}}%
+g_{jb}(x)\Gamma _{ra}^{b}(x)\right) \right] _{[rj]}
\end{equation}%
as the first component of $\mathcal{R}$ and rewrite (67) in the form

\begin{equation}
I_{1}(\mathbf{g};y)-(x,f(x),f_{1}(x))_{\ast }I_{1}(\mathbf{g};x)=0
\end{equation}

As a very crucial fact, if $\Gamma =(\Gamma _{jk}^{i})$ are the Christoffel
symbols of the Levi-Civita connection $\nabla ,$ then

\begin{equation}
\nabla _{r}g_{jk}=\frac{\partial g_{jk}}{\partial x^{r}}+g_{ja}(x)\Gamma
_{rk}^{a}(x)+g_{ka}(x)\Gamma _{rj}^{a}(x)
\end{equation}%
and using (68) and (70) we easily check that

\begin{equation}
I_{rj}^{k}(\mathbf{g};x)=\left[ g^{ak}\nabla _{r}g_{ja}\right] _{[rj]}
\end{equation}

\textit{Therefore }$I_{1}(\mathbf{g};x)=0$\textit{\ if }$\Gamma =(\Gamma
_{jk}^{i})$\textit{\ are the Christoffel symbols, i.e., this assumption
makes the first component of }$\mathcal{R}$ \textit{vanish!!}

To deduce the second component of $\mathcal{R},$ we differentiate the second
equation of (45) and substitute back from it... and what we get is clearly
(47)+(48)+(49) which we will rewrite here:

\begin{eqnarray}
\mathcal{R}_{rj,k}^{i}(x,f(x),f_{1}(x))\overset{def}{=} &&\left[ \frac{%
\partial \Gamma _{jk}^{i}(f(x))}{\partial y^{r}}+\Gamma
_{jb}^{i}(f(x))g_{rk}^{b}(f(x))\right] _{[rj]} \\
&&-g_{j}^{a}g_{r}^{c}g_{k}^{d}\left[ \frac{\partial \Gamma _{ad}^{e}(x)}{%
\partial x^{c}}+\Gamma _{ab}^{e}(x)\Gamma _{cd}^{b}(x)\right]
_{[ca]}f_{e}^{i}(x)=0  \notag
\end{eqnarray}

\begin{equation}
I_{2}(\mathbf{g};x)\overset{def}{=}I_{rj,k}^{i}(\mathbf{g};x)\overset{def}{=}%
\left[ \frac{\partial \Gamma _{jk}^{i}(x)}{\partial x^{r}}+\Gamma
_{ja}^{i}(x)\Gamma _{rk}^{a}(x)\right] _{[rj]}
\end{equation}

\begin{equation}
\mathcal{R(}x,y,f_{1})=I_{2}(\mathbf{g};y)-(x,y,f_{1})_{\ast }I_{2}(\mathbf{g%
};x))=0
\end{equation}

Now (45) is uniquely locally integrable if and only if

\begin{eqnarray}
I_{1}(\mathbf{g};y)-(x,f(x),f_{1}(x))_{\ast }I_{1}(\mathbf{g};x)\overset{def}%
{=}\mathcal{R}_{1}(x,f(x),f_{1}(x)) &=&0 \\
I_{2}(\mathbf{g};y)-(x,f(x),f_{1}(x))_{\ast }I_{2}(\mathbf{g};x)\overset{def}%
{=}\mathcal{R}_{2}(x,f(x),f_{1}(x)) &=&0  \notag
\end{eqnarray}

We define

\begin{equation}
\mathcal{R}\overset{def}{=}(\mathcal{R}_{1},\mathcal{R}_{2})
\end{equation}%
and 
\begin{equation}
I(\mathbf{g};x)\overset{def}{=}(I_{1}(\mathbf{g)},I_{2}(\mathbf{g}%
))=(I_{rj}^{i}(\mathbf{g};x),I_{rj,k}^{i}(\mathbf{g};x))
\end{equation}

To clarify the meaning of $I(\mathbf{g};x),$ a section of $%
J_{1}(T)\rightarrow M$ is of the form $(\xi ^{i},\xi _{k}^{i})$ and using
regular coordinates at $x,$ it is not difficult to check that $I(\mathbf{g}%
;x)$ is a $2$-form at $x$ with values in the fiber of $J_{1}(T)\rightarrow M$
over $x$ (\textit{but not necessarily in the fiber of }$K_{1}(\mathbf{g}%
)\rightarrow M$ over $x!!)$ i.e., $I(\mathbf{g};x)$ is a section of the
vector bundle $\wedge ^{2}(T^{\ast })\otimes J_{1}(T)\longrightarrow M.$

Now we rewrite (75) in compact form as

\begin{equation}
\mathcal{R}(x,y,f_{1})\overset{def}{=}I(\mathbf{g};y)-(x,y,f_{1})_{\ast }I(%
\mathbf{g};x)
\end{equation}%
and it remains to clarify the meaning of $(x,y,f_{1})_{\ast }$ in (78). We
recall that the vector bundle $J_{1}(T)\rightarrow M$ is associated with $%
\mathcal{U}_{2}$ in the sense that a bisection $(x,f(x),f_{1}(x),f_{2}(x))$
of $\mathcal{U}_{2}$ induces an isomorphism between the fibers $%
(x,f(x),f_{1}(x),f_{2}(x))_{\ast }:J_{1}(T)^{x}\rightarrow J_{1}(T)^{y}.$
Now the bisection $(x,f(x),f_{1}(x))$ of $\mathcal{K}_{1}(\mathbf{g})$ lifts
by (66) to a bisection of $\mathcal{K}_{2}(\mathbf{g})$ and induces the
isomorphism $(x,y,f_{1})_{\ast }:J_{1}(T)^{x}\rightarrow J_{1}(T)^{y}$ which
restricts to $(x,y,f_{1})_{\ast }:K_{1}(\mathbf{g})\rightarrow J_{1}(T)^{y}.$
Therefore we obtain the map

\begin{equation}
\mathcal{R}:\mathcal{K}_{1}(\mathbf{g})\longrightarrow \wedge ^{2}(T^{\ast
})\otimes J_{1}(T)
\end{equation}%
which gives the full integrability conditions of (58)+(59). We emphasize
again that $\mathcal{R}$ does not necessarily take values in the smaller
bundle $\wedge ^{2}(T^{\ast })\otimes K_{1}(\mathbf{g})$ (therefore our last
statement in the second paragraph on page 192 of [Or1] is incorrect). We
observe the following remarkable fact: $I(\Gamma )$ is a tensor in affine
geometry and has an invariant meaning whereas $I_{2}(\mathbf{g)}$ is a
component of the second order object $\mathcal{R}$ in coordinates and has no
invariant meaning alone unless $I_{1}(\mathbf{g})=0!!$

Now the full integrability conditions of (63)+(64) are obtained by the
linearization of (58)+(59) and is given by

\begin{equation}
\mathfrak{R}:K_{1}(\mathbf{g})\longrightarrow \wedge ^{2}(T^{\ast })\otimes
J_{1}(T)
\end{equation}%
where $\mathfrak{R}$ is a section of the vector bundle $\wedge ^{2}(T^{\ast
})\otimes Hom(K_{1}(\mathbf{g}),J_{1}(T)).$

Therefore we state

\begin{proposition}
(Lie's 3'rd Theorem) The following are equivalent

1) $\mathcal{R}=0$ on $\mathcal{K}_{1}(\mathbf{g})$

2) $\mathcal{K}_{2}(\mathbf{g})$ is uniquely locally integrable

3) $\mathcal{K}_{2}(\mathbf{g}),$ which leaves $\mathbf{g}$ invariant by
definition, leaves also $I(\mathbf{g})$ invariant

4) $\mathfrak{R}=0$

5) $K_{2}(\mathbf{g})$ is uniquely locally integrable

6) $\mathcal{L}_{\varepsilon \xi _{1}}I(\mathbf{g})=0$ for all sections $\xi
_{1}$ of $K_{1}(\mathbf{g})\rightarrow M$
\end{proposition}

By Proposition 21, we know that an affine PHG is flat if and only if the
torsionfree affine connection $\Gamma $ is flat. The answer to the following
question will point at a new phenomenon for a Riemannian PHG.

\textbf{Q : }What is the meaning of the conditions of Proposition 21 in
terms of the metric $g=(g_{ij})?$

The very neat answer is given by

\begin{proposition}
For a Riemannian PHG $\mathcal{K}_{2}(\mathbf{g}),$ the following are
equivalent

1) One of the conditions of Proposition 21 holds

2) $\Gamma $ is the Levi-Civita connection and the metric $g$ has constant
curvature.
\end{proposition}

To prove Proposition 22, we first claim $\mathcal{R}=0$ implies that $\Gamma
=(\Gamma _{jk}^{i})$ are the Christoffel symbols. Let $y=f(x)$ be local
diffeomorphism satisfying (58)+(59) with the initial conditions (60)+(61).
Differentiating (58) with respect to $x^{r}$ and evaluating at $x=\overline{x%
},$ we get

\begin{equation}
\frac{\partial g_{ab}(\overline{y})}{\partial y^{c}}\overline{f}_{r}^{c}%
\overline{f}_{j}^{a}\overline{f}_{k}^{b}+g_{ab}(\overline{y})\overline{f}%
_{rj}^{a}\overline{f}_{k}^{b}+g_{ab}(\overline{y})\overline{f}_{j}^{b}%
\overline{f}_{rk}^{a}=\frac{\partial g_{jk}(\overline{x})}{\partial x^{r}}
\end{equation}

Shifting the indices $r,j,k$ in (81) using the Gauss trick and defining

\begin{eqnarray}
&&\widetilde{\Gamma }_{rj}^{k}(\overline{x})\overset{def}{=}\frac{1}{2}%
g^{ka}\left( \frac{\partial g_{ja}(\overline{x})}{\partial x^{r}}-\frac{%
\partial g_{rj}(\overline{x})}{\partial x^{a}}+\frac{\partial g_{ar}(%
\overline{x})}{\partial x^{j}}\right) \\
&=&\text{ the Christoffel symbols of }g  \notag
\end{eqnarray}%
(81) becomes after some computation

\begin{equation}
\widetilde{\Gamma }_{ab}^{i}(\overline{y})\overline{f}_{j}^{a}\overline{f}%
_{k}^{b}=\widetilde{\Gamma }_{jk}^{a}(\overline{x})\overline{f}_{a}^{i}+%
\overline{f}_{jk}^{i}
\end{equation}

Therefore (60)+(61) and (60)+(83) have the same solutions for \textit{all }%
initial conditions (60)+(61). Subtracting (83) from (61), setting $\overline{%
y}=\overline{x}$ and using regular coordinates around $\overline{x},$ we
conclude that $O(n)$ fixes the $(1,2)$-tensor $\Gamma _{jk}^{i}(\overline{x}%
)-\widetilde{\Gamma }_{jk}^{i}(\overline{x}).$ Since $1+2=3$ is odd, we
conclude $\Gamma _{jk}^{i}(\overline{x})-\widetilde{\Gamma }_{jk}^{i}(%
\overline{x})=0$ as claimed (This derivation actually uses a much weaker
assumption than $\mathcal{R}=0,$ see the next section). Therefore $I_{1}(%
\mathbf{g})=0$ by (71) and we deduce from (73) that

\begin{equation}
I_{2}(\mathbf{g})=R=\text{ the Riemann curvature tensor of }g
\end{equation}

Now it follows from Proposition 21 that $R_{kj,m}^{l}$ and therefore $%
R_{kj,lm}$ is fixed pointwise by $O(n)$ and therefore can be expressed in
terms of the components of $g=(g_{ij})$ since all tensor invariants of $O(n)$
are of this form. We define

\begin{equation}
\overline{R}_{kj,lm}\overset{def}{=}g_{lk}g_{jm}-g_{lj}g_{km}
\end{equation}%
and check the curvature identities

\begin{eqnarray}
\overline{R}_{kj,lm} &=&-\overline{R}_{jk,lm}  \notag \\
\overline{R}_{kj,lm} &=&-\overline{R}_{kj,ml} \\
\overline{R}_{kj,lm}+\overline{R}_{lk,jm}+\overline{R}_{jl,km} &=&0  \notag
\end{eqnarray}

Since the vector space of the $O(n)$-invariant tensors $\overline{R}_{kj,lm}$
satisfying (86) is spanned by $\overline{R}_{kj,lm},$ we conclude

\begin{equation}
I_{kj,lm}(g,\Gamma ;\overline{x})=R_{kj,lm}(\overline{x})=c(\overline{x})%
\overline{R}_{kj,lm}=c(\overline{x})(g_{lk}g_{jm}-g_{lj}g_{km})
\end{equation}%
for some scalar $c(\overline{x})$. Since $\overline{R}$ and $R$ are both $%
\mathcal{K}_{1}(\mathbf{g})$-invariant on $M,$ it follows that $c(\overline{x%
})$ does not depend on $\overline{x}$ and (87) becomes the constant
curvature condition. Hence 1) implies 2). Conversely, if $\Gamma =\widetilde{%
\Gamma },$ then $I_{1}(\mathbf{g})=0$ by (71) and the above argument showes
that 1) follows from 2), finishing the proof of Proposition 15.

We note here that the equivalence of 4) of Proposition 21 and 2) of
Proposition 22 is shown also in [Po], 254-255 making heavy use of Spencer
cohomology.

\section{Connections}

The equations (81) are obtained differentating (60) and substituting jet
variables. The equations (60)+(81) are called the (first) prolongation of
(60) and denoted by $pr(\mathcal{K}_{1}(\mathbf{g})).$ The proof of the
first part of Proposition 22 shows actually

\begin{eqnarray}
\mathcal{K}_{2}(\mathbf{g}) &\subset &pr(\mathcal{K}_{1}(\mathbf{g}))\iff 
\mathcal{K}_{2}(\mathbf{g})=pr(\mathcal{K}_{1}(\mathbf{g})) \\
&\iff &\Gamma =\text{ Christoffel symbols}  \notag
\end{eqnarray}

$\mathcal{K}_{2}(\mathbf{g})$ is called $1$-flat (called $1$-torsionfree in
[Or1], unfortunately another misleading terminology) if (88) is satisfied.
Of course, $\mathcal{R}=0$ implies $1$-flatness of $\mathcal{K}_{2}(\mathbf{g%
}).$ For an affine PHG $\mathcal{H}_{2}(\Gamma ),$ we have

\begin{equation}
\mathcal{H}_{2}(\Gamma )\subset pr(\mathcal{H}_{1}(\Gamma ))=pr(\mathcal{U}%
_{1})=\mathcal{U}_{2}
\end{equation}%
and therefore $\mathcal{H}_{2}(\Gamma )$ is always $1$-flat. Similarly, for
an absolute parallelism $\mathcal{N}_{1}(w),$ we have

\begin{equation}
\mathcal{N}_{1}(w)\subset pr(\mathcal{N}_{0}(w))=pr(M\times M)=\mathcal{U}%
_{1}
\end{equation}%
and $\mathcal{N}_{1}(w)$ is always $0$-flat. Even though $1$-flatness of $%
\mathcal{K}_{2}(\mathbf{g})$ is equivant to the condition $\Gamma $ =
Christoffel symbols (the reason why (88) is called $1$-torsionfreeness in
[Or1] !!), we observe \textit{that }$0$\textit{-flatness of }$\mathcal{N}%
_{1}(w)$\textit{\ has nothing to do with the torsions of the affine
connections }$\nabla $\textit{\ and }$\widetilde{\nabla }!!$\textit{\
Similarly, }$1$\textit{-flatness of }$\mathcal{H}_{2}(\Gamma )$\textit{\ has
nothing to do with the torsionfreeness of the affine connection }$\Gamma .$

$i$-flatness gurantees that the first operators in Spencer sequences take
values in the "right \ spaces" as follows: Starting with $\mathcal{N}%
_{1}(w), $ we have the first nonlinear Spencer sequence

\begin{equation}
\mathcal{N}_{0}(w)=M\times M\overset{\mathcal{D}_{1}}{\longrightarrow }%
\wedge ^{1}(T^{\ast })\otimes T\overset{\mathcal{D}_{2}}{\longrightarrow }%
\wedge ^{2}(T^{\ast })\otimes T
\end{equation}%
and the first three terms of the linear Spencer sequence

\begin{equation}
T\overset{D_{1}}{\longrightarrow }\wedge ^{1}(T^{\ast })\otimes T\overset{%
D_{2}}{\longrightarrow }\wedge ^{2}(T^{\ast })\otimes T\longrightarrow ....
\end{equation}

Coordinate description of the operators in (91) are given in [Or1] and $%
D_{1}=\nabla $ defined by (41). \textit{Note that }$D_{1}\neq \widetilde{%
\nabla }=\mathcal{L}!!$ Now $\mathcal{R}$ and $\mathfrak{R}$ drop out of the
compositions $\mathcal{D}_{2}\circ \mathcal{D}_{1}$ and $D_{2}\circ D_{1}$
respectively and Proposion 21 asserts that (91) is locally exact if and only
if (92) is locally exact. From our standpoint, the reason why $D_{1}$ turns
out to be a connection is that $\mathcal{N}_{1}(w)$ is $0$-flat as will
become clear below. Similar remarks apply to an affine PHG $\mathcal{H}%
_{2}(\Gamma )$ as follows: We have the nonlinear and linear Spencer sequences

\begin{equation}
\mathcal{H}_{1}(\Gamma )=\mathcal{U}_{1}\overset{\mathcal{D}_{1}}{%
\longrightarrow }\wedge ^{1}(T^{\ast })\otimes H_{1}(\Gamma )\overset{%
\mathcal{D}_{2}}{\longrightarrow }\wedge ^{2}(T^{\ast })\otimes H_{1}(\Gamma
)
\end{equation}

\begin{equation}
H_{1}(\Gamma )\overset{D_{1}}{\longrightarrow }\wedge ^{1}(T^{\ast })\otimes
H_{1}(\Gamma )\overset{D_{2}}{\longrightarrow }\wedge ^{2}(T^{\ast })\otimes
H_{1}(\Gamma )\longrightarrow ....
\end{equation}%
and $\mathfrak{R}$ is the curvature of the connection $D_{1}$ and is a
section of $\wedge ^{2}(T^{\ast })\otimes Hom(H_{1}(\Gamma ),H_{1}(\Gamma
))= $ $\wedge ^{2}(T^{\ast })\otimes Hom(J_{1}(T),J_{1}(T)).$ Note again
that $\mathfrak{R}$ is by no means the curvature of the torsionfree affine
connection $\Gamma $ defined on the tangent bundle $T\rightarrow M.$

Something very interesting happens for the Riemannian PHG $\mathcal{K}_{2}(%
\mathbf{g}).$ Lifting a section of $\mathcal{K}_{1}(\mathbf{g})$ to a
section of $\mathcal{K}_{2}(\mathbf{g})$ and then mapping by the Spencer
operator $\mathcal{D}_{1}$ (see the Appendix of [Or1] and [P1], [P2] for the
definition of $\mathcal{D}_{1}),$ we get the map

\begin{equation}
\mathcal{K}_{1}(\mathbf{g})\overset{\mathcal{D}_{1}}{\longrightarrow }\wedge
^{1}(T^{\ast })\otimes J_{1}(T)
\end{equation}%
but \textit{not necessarily a map}

\begin{equation}
\mathcal{K}_{1}(\mathbf{g})\overset{\mathcal{D}_{1}}{\longrightarrow }\wedge
^{1}(T^{\ast })\otimes K_{1}(\mathbf{g})\subset \wedge ^{1}(T^{\ast
})\otimes J_{1}(T)
\end{equation}

Linearizations of (95) and (96) are given by

\begin{equation}
K_{1}(\mathbf{g})\overset{D_{1}}{\longrightarrow }\wedge ^{1}(T^{\ast
})\otimes J_{1}(T)
\end{equation}

\begin{equation}
K_{1}(\mathbf{g})\overset{D_{1}}{\longrightarrow }\wedge ^{1}(T^{\ast
})\otimes K_{1}(\mathbf{g})
\end{equation}

We observe that $D_{1}$ in (97) is \textit{not a connection on }$K_{1}(%
\mathbf{g})!!$\textit{\ }Now it is easy to check that (88) holds if and only
if (98) holds. With this assumption, $D_{1}$ becomes a connection and we get
the first three terms of the linear Spencer sequence

\begin{equation}
K_{1}(\mathbf{g})\overset{D_{1}}{\longrightarrow }\wedge ^{1}(T^{\ast
})\otimes K_{1}(\mathbf{g})\overset{D_{2}}{\longrightarrow }\wedge
^{2}(T^{\ast })\otimes K_{1}(\mathbf{g})\longrightarrow ....
\end{equation}

Now $\mathfrak{R}$ becomes the curvature of the connection $D_{1}$ and is a
section of $\wedge ^{2}(T^{\ast })\otimes Hom(K_{1}(\mathbf{g}),K_{1}(%
\mathbf{g})).$ \textit{However, recall that (88) forces }$I_{1}(\mathbf{g}%
)=0 $ which makes $\mathfrak{R}$\textit{\ a section of }$\wedge ^{2}(T^{\ast
})\otimes Hom(K_{1}(g),T^{\ast }\otimes T)!$

We will conclude with three remarks.

1) Since $\mathfrak{R}$ is a section of $\wedge ^{2}(T^{\ast })\otimes
Hom(K_{1}(\mathbf{g}),J_{1}(T))$ in general, one may object that $\mathfrak{R%
}$ is not intrinsic to $\mathcal{K}_{2}(\mathbf{g})$ and (88) is necessary
to make it intrinsic. However, $\mathcal{K}_{2}(\mathbf{g})$ is \textit{by
definition }a subgoupoid of $\mathcal{U}_{2}$ and it is not possible to
seperate it from $\mathcal{U}_{2}$ and define it as an abstract structure,
i.e., $\mathfrak{R}$ is intrinsic to jets. Indeed, an abstract $G$-principal
bundle is not always a $G$-structure defined as a reduction of the principal
bundle $\mathcal{U}_{2}^{e,\bullet }\rightarrow M$ (like $\mathcal{K}_{2}(%
\mathbf{g})).$ For instance, the concept of torsion emerges from connections
on $G$-reductions and makes no sense (unless we introduce further structure)
for connections on abstract $G$-principal bundles.

2) The above linear analysis becomes more intriguing (even catastrophic) in
the nonlinear case when we attempt to interpret $\mathcal{R}$ as the
curvature of a torsionfree connection on the principal bundles $\mathcal{U}%
_{1}^{e,\bullet }\rightarrow M$ (affine) and $\mathcal{K}_{1}(\mathbf{g}%
)^{e,\bullet }\rightarrow M$ (riemannian) and will be omitted here.

3) We recall that $\mathcal{L}_{\varepsilon \xi _{k}}$ is not a connection
for $k\geq 1.$ This operator is used to localize the global sequences (92),
(94), (99) and more generally the linear Spencer sequences arising from
PHG's. This localization generalizes the well known process of passing from
de Rham cohomology to Lie algebra cohomology by localizing the forms in the
sequence. The details of this construction for absolute parallelism are
given in Chapter 11 of [Or1].

\section{Reduction to parallelism}

We start with an affine PHG $\mathcal{H}_{2}(\Gamma ).$ We fix a basepoint $%
e\in M$ and consider the right principal bundle $\pi :\mathcal{H}_{1}(\Gamma
)^{e,\bullet }=$ $\mathcal{U}_{1}^{e,\bullet }\rightarrow M$ with the
structure group $\mathcal{U}_{1}^{e,e}\cong GL(n,\mathbb{R)}$ whose fiber $%
\pi ^{-1}(x)$ over $x$ is the set $\mathcal{U}_{1}^{e,x}$ of all $1$-arrows
from $e$ to $x.$ Let $\widetilde{x}\in \mathcal{U}_{1}^{e,x},$ $\pi (%
\widetilde{x})=x$ and consider the tangent space $T_{\widetilde{x}}(\mathcal{%
U}_{1}^{e,\bullet })$ of $\mathcal{U}_{1}^{e,\bullet }$ at $\widetilde{x}$
which is easily seen to be canonically isomorphic to the fiber $J_{1}(T)^{x}$
of $J_{1}(T)\rightarrow M$ over $x,$ i.e.,

\begin{equation}
T_{\widetilde{x}}(\mathcal{U}_{1}^{e,\bullet })\cong J_{1}(T)^{x}\text{ \ \
\ }\pi (\widetilde{x})=x
\end{equation}%
The isomorphism in (100) is obtained by lifting the $1$-parameter subgroup
defined by $\xi _{1}\in J_{1}(T)^{x}$ to $\mathcal{U}_{1}^{e,\bullet }$ by
composition \textit{at the target} and differentiating at $t=0.$ In
coordinates, it is given by

\begin{equation}
\xi ^{a}(x)\frac{\partial }{\partial x^{a}}+\xi _{b}^{a}(x)f_{c}^{b}(x)\frac{%
\partial }{\partial f_{c}^{a}}\longleftrightarrow (\xi ^{i}(x),\xi
_{j}^{i}(x))
\end{equation}%
where $\widetilde{x}=(e^{i},x^{i},f_{j}^{i})\in \mathcal{U}_{1}^{e,x}.$ Note
that the action of $\mathcal{U}_{1}^{e,e}\cong GL(n,\mathbb{R)}$ on $%
\mathcal{U}_{1}^{e,\bullet }$ by composition \textit{at the source }commutes
with the above isomorphism. Therefore, some $\xi _{1}\in J_{1}(T)^{x}$
defines tangent vectors at all points in the fiber $\pi ^{-1}(x)$ which are
invariant by the action of $\mathcal{U}_{1}^{e,e}.$

Now let $\widetilde{x},$ $\widetilde{y}\in \mathcal{U}_{1}^{e,\bullet },$ $%
\pi (\widetilde{x})=x,$ $\pi (\widetilde{y})=y$ and consider the $1$-arrow $%
\widetilde{y}\circ \widetilde{x}^{-1}$ on $M$ from $x$ to $y.$ Therefore $%
\varepsilon (\widetilde{y}\circ \widetilde{x}^{-1})$ is a $2$-arrow of $%
\mathcal{H}_{2}(\Gamma )\cong \mathcal{H}_{1}(\Gamma )=\mathcal{U}_{1}$ on $%
M $ from $x$ to $y.$ Recalling that $H_{1}(T)=J_{1}(T)$ is associated with $%
\mathcal{U}_{2},$ the $2$-arrow $\varepsilon (\widetilde{y}\circ \widetilde{x%
}^{-1})$ induces an isomorphism

\begin{equation}
\varepsilon (\widetilde{y}\circ \widetilde{x}^{-1}):J_{1}(T)^{x}%
\longrightarrow J_{1}(T)^{y}
\end{equation}

In view of (100), (102) becomes

\begin{equation}
\varepsilon (\widetilde{y}\circ \widetilde{x}^{-1}):T_{\widetilde{x}}(%
\mathcal{U}_{1}^{e,\bullet })\longrightarrow T_{\widetilde{y}}(\mathcal{U}%
_{1}^{e,\bullet })
\end{equation}

We conclude from (103) that $\varepsilon (\widetilde{y}\circ \widetilde{x}%
^{-1}),$ which is a $2$-arrow on $M$ from $x$ to $y,$ is at the same time a $%
1$-arrow on $\mathcal{U}_{1}^{e,\bullet }$ from $\widetilde{x}$ to $%
\widetilde{y}$ $!!$ We easily check that these $1$-arrows are closed under
composition and inversion and we conclude

\begin{proposition}
The splitting (44) defines an absolute parallelism on $\mathcal{H}%
_{1}^{e,\bullet }(\Gamma )=\mathcal{U}_{1}^{e,\bullet }.$
\end{proposition}

\bigskip Using the notation of Section 1, we will write the absolute
parallelism in Proposition (23) in the form

\begin{equation}
\varepsilon _{\alpha }^{\beta }(\widetilde{x},\widetilde{y})
\end{equation}%
where the indices $\alpha ,\beta $ refer to the components of the tangent
vectors at $\widetilde{x},\widetilde{y}.$ It follows from (44) and (103)
that $\varepsilon _{\alpha }^{\beta }(\widetilde{x},\widetilde{y})$ depends
on $\Gamma (x),$ $\Gamma (y)$ but we will not bother here to write down the
explicit form of (104). Now (104) explains the reason for the strong analogy
between our computations in Sections 2, 3 but much more importantly, it
allows us to carry everything done for absolute parallelism in Part 2 of
[Or1] over to the affine case (like Chern-Simons classes, homogeneous
flow...etc).

When we attempt to generalize (104) to a Riemannian PHG, a problem arises:
Now $\mathcal{K}_{1}(\mathbf{g})^{e,\bullet }\varsubsetneq \mathcal{U}%
_{1}^{e,\bullet }$ and $J_{1}(T)\rightarrow M,$ which is associated with $%
\mathcal{U}_{2},$ is clearly associated also with $\mathcal{K}_{2}(\mathbf{g}%
)\subset \mathcal{U}_{2},$ i.e., a $2$-arrow $\varepsilon (f^{x,y})$ of $%
\mathcal{K}_{2}(\mathbf{g})$ from $x$ to $y$ induces an isomorphism $%
\varepsilon (f^{x,y}):J_{1}(T)^{x}\rightarrow $ $J_{1}(T)^{y}.$ However,
this isomorphism may not restrict to an isomorphism $K_{1}(\mathbf{g}%
)^{x}\rightarrow K_{1}(\mathbf{g})^{y},$ i.e., $K_{1}(\mathbf{g})\rightarrow
M$ need not be associated with $\mathcal{K}_{2}(\mathbf{g}):$ All we can
conclude is the injection $\varepsilon (f^{x,y}):K_{1}(\mathbf{g}%
)^{x}\rightarrow $ $J_{1}(T)^{y}.$ In fact, using the definitions, it is not
difficult to show the following

\begin{proposition}
For a Riemannian PHG the following are equivalent

1) $K_{1}(\mathbf{g})\rightarrow M$ is associated with $\mathcal{K}_{2}(%
\mathbf{g})$

2) $\mathcal{K}_{2}(\mathbf{g})$ is $1$-flat, i.e., $\mathcal{K}_{2}(\mathbf{%
g})\subset prK_{1}(\mathbf{g})$ or equivalently $\Gamma =(\Gamma _{jk}^{i})$
are the Christoffel symbols.
\end{proposition}

Fortunately, there is an easy way out of this difficulty.\ 

\begin{definition}
Let $L\subset S$ be a submanifold. Then $L$ is parallelizable relative to $S$
if for any $x,y\in L,$ there is a unique $1$-arrow $\varepsilon (x,y)$ of $S$
(not necessarily of $L!)$ such that these $1$-arrows of $S$ are closed under
composition and inversion.
\end{definition}

We call $\varepsilon $ an $S$-splitting of $L.$ If the $1$-arrows $%
\varepsilon (x,y)$ of $S$ restrict to the $1$-arrows of $L,$ i.e., if the
maps $\varepsilon (x,y):T(S)^{x}\rightarrow T(S)^{y}$ restrict to $%
\varepsilon (x,y):T(L)^{x}\rightarrow T(L)^{y}$ (our fiber notation $A^{x}$
forces us to denote the tangent space $T_{x}(S)$ by $T(S)^{x}),$ then the $S$%
-splitting $\varepsilon $ on $L$ becomes a true splitting and $L$ becomes
absolutely parallelizable.

Now the key fact is that if $L$ is parallelizable relative to $S,$ then both
sides of (14) are well defined and are contained in the tangent space $%
T(S)^{f(x)}.$ Therefore, the equality in (14) makes perfect sense with the
only difference that the index $i$ on the LHS of (14) refers to the tangent
space of $L$ whereas on the RHS of (14) it refers to the tangent space of $%
S. $ Now all the computations of Section 2 work through if we replace
absolute parallelism with relative parallelism. For instance, the linear
curvature $\mathfrak{R}_{jk,l}^{(i)}$ becomes a $2$-form on $L$ with values
in $Hom(T(L),T(S))...etc.$

Our standard example of relative parallelism is $L=\mathcal{K}_{1}(\mathbf{g}%
)^{e,\bullet }\subset \mathcal{U}_{1}^{e,\bullet }=S$ and $\varepsilon $ is
as in (104) which we may write now as $\varepsilon _{\alpha }^{(\beta )}(%
\widetilde{x},\widetilde{y}).$ Note that our choice of $S$ is \textit{%
canonical. }As in the affine case, we can now reduce the study of a
Riemannian PHG to that of relative parallelism and furthermore carry all
constructions of Part 2 of [Or1] over to the Riemannian case.

\section{\protect\bigskip General PHG's}

When we attempt to define and study a general PHG, the first problem we face
is the construction of the structure object. Once this is done, the rest
follows as in the case of parallelizable, affine and Riemannian PHG's along
the same lines.

Let $G$ be a Lie group acting transitively on a smooth manifold $M.$ Fixing
a point $e\in M,$ we can identify this action with the (say) left action of $%
G$ on $G/H$ where $H$ is the stabilizer at $e.$ We fix some $p,q\in M,$ $%
g\in G$ with $g(p)=q$ and define

\begin{equation}
G_{k}(p,q;g)\overset{def}{=}\{f\in G\mid f(p)=q,\text{ }%
j_{k}(f)^{p,q}=j_{k}(g)^{p,q}\}
\end{equation}

Obviously $g\in G_{k}(p,q;g)$ for all $k\geq 0$ and $G_{k+1}(p,q;g)\subset
G_{k}(p,q;g).$ With some mild assumptions (like $G$ is connected and acts
effectively) we can show (see [Or1]) the existence of a smallest integer $m$
with the property that $G_{m+1}(p,q;g)=\{g\}$ and furthermore $m$ is \textit{%
independent }of $p,q$ and $g.$ In particular, we may choose $p=q$
arbitrarily. In short, any transformation of $G$ is \textit{globally
determined on }$M$ by \textit{any }of its $m$-arrows. It follows that above
any $m$-arrow, there is a unique $(m+1)$-arrow. Indeed, since $g\in G$ is
determined by $j_{m}(g)^{x,g(x)}$ for any $x\in M,$ $j_{m+1}(g)^{x,g(x)}%
\overset{def}{=}$ $\varepsilon (j_{m}(g)^{x,g(x)})$ is the unique $(m+1)$%
-arrow above $j_{m}(g)^{x,g(x)}.$ In this way, we obtain a transitive
subgroupoid $\mathcal{P}_{m+1}\subset $ $\mathcal{U}_{m+1}$ with the
property $\mathcal{P}_{m+1}=\varepsilon (\mathcal{P}_{m})\cong \mathcal{P}%
_{m}.$ It turns out that $\mathcal{P}_{m+1}$ defines a first order nonlinear
system of PDE's on the universal pseudogroup $Diff_{loc}(M)$ of all local
diffeomorphisms of $M$ whose unique solutions are the restrictions of
actions of the elements of $G.$ Therefore, any $(m+1)$-arrow of $\mathcal{P}%
_{m+1}$ integrates uniquely to a global transformation of $G.$ The PHG $%
\mathcal{P}_{m+1}$ is a flat model with $\mathcal{R}=0.$ The idea is now to
forget the solution space $G$ and define $\mathcal{P}_{m+1}$ as an
independent structure on a smooth $N$ with $\dim N=\dim M$ with the "same
stabilizer $H"$ and with curvature $\mathcal{R}$ which will be the
obstruction to the existence of (unique) local solutions. In order to do
this, we need first to construct the structure object of our flat model $%
\mathcal{P}_{m+1}$ on $M.$

The above construction of the integer $m$ gives now the injective map

\begin{equation}
j_{m+1}:H\longrightarrow j_{m+1}(h)^{e,e}\in \mathcal{U}_{m+1}^{e,e}\cong
G_{m+1}(n)
\end{equation}%
for $h\in H$ and $j_{m}(H)\cong j_{m+1}(H)=\varepsilon (j_{m}(H)).$ It is
easy to construct examples of action pairs $(G,M)$ with arbitrarily large $m$
using graded Lie algebras of vector fields (see [D] and the references
therein) but it seems to us that the structure of such pairs is far from
being well understood. Clearly, $m=0$ for parallelism, $m=1$ for affine and
Riemannian PHG's and $m=2$ for projective and conformal PHG's.

Now suppose that $G$ is an algebraic group and $H\subset G$ an algebraic
subgroup. Recalling that $G_{m+1}(n)$ is an affine algebraic group, the
first problem is the following question

\textbf{Q1: }Is $j_{m+1}(H)=\varepsilon (j_{m}(H))\subset G_{m+1}(n)$ an
algebraic subgroup?

Now consider the left coset space $G_{m+1}(n)/j_{m+1}(H)$ and the action of $%
G_{m+1}(n)$ on $G_{m+1}(n)/j_{m+1}(H).$ Assuming that the answer to \textbf{%
Q1 }is affirmative, we now ask

\textbf{Q2: }Does there exist polynomial functions on $G_{m+1}(n)$
seperating the left cosets of $j_{m+1}(H),$ i.e., a polynomial injective map 
\begin{equation}
\Omega :G_{m+1}(n)/j_{m+1}(H)\longrightarrow \mathbb{R}^{s}
\end{equation}%
for some $s.$

The answers to \textbf{Q1, Q2 }are affirmative for parallelizable, affine
and Riemannian PHG's. In fact, the components of the structure objects $%
w=(w_{j}^{i}),$ $\Gamma =(\Gamma _{jk}^{i}),$ $\mathbf{g}=(g_{ij},\Gamma
_{jk}^{i})$ give the required imbeding (107) and $s$ is minimal in all these
cases (see [Or1] for more details). We believe that the answers to \textbf{%
Q1, Q2 }are affirmative if $H\subset G$ are algebraic groups and therefore $%
G_{m+1}(n)/j_{m+1}(H)$ is always an affine variety.

Now assuming we found the polynomial function $\Omega =(\Omega ^{\alpha })$
whose components parametrize the left coset space $G_{m+1}(n)/j_{m+1}(H),$
the action of $G_{m+1}(n)$ on $G_{m+1}(n)/j_{m+1}(H)$ gives the
transformation rule of the components $(\Omega ^{\alpha }).$ Now it is easy
to check that the $(m+1)$-arrows of our flat $\mathcal{P}_{m+1}$ on $M$ are
those $(m+1)$-arrows in $\mathcal{U}_{m+1}$ that preserve the geometric
object $\Omega .$ Since $\mathcal{P}_{m+1}$ is flat, $\Omega $ is subject to
some "integrability conditions". Equivalently, the transformations of $G$
preserve both $\Omega $ and its integrability object $I(\Omega )$ and as a
remarkable fact, $\Omega $ and $I(\Omega )$ now drop out of the
all-important recursion formulas in the Fels-Olver theory of moving frames
([FO2]) that we will briefly mention in the next section.

Having the geometric object $\Omega =(\Omega ^{\alpha })$ at our disposal,
we now start anew on a smooth manifold $N$ with $\dim N=\dim M$ by
postulating the existence of some $\Omega =(\mathcal{\Omega }^{\alpha }(x))$
on $N$ whose components are subject to the above transformation rule and
define $\mathcal{S}_{m+1}\subset \mathcal{U}_{m+1}$ as the invariance
subgroupoid of $\Omega .$ It is now a straightforward matter to linearize $%
\mathcal{S}_{m+1}$ and prove Lie's 3'rd Theorem. Many highly nontrivial
questions arise, but it seems pointless at this stage to eloborate further
on the theory before making a detailed study of projective and conformal
structures and understanding the advantages and disadvantages of the present
theory of geometric structures.

\section{Moving frames}

In [FO1], [FO2] Olver and Fels introduced a new theoretical foundation of
the moving frame method most closely associated with Elie Cartan which is
very simple to apply and amazingly powerful. The wide range of new
applications of this new approach (see the survey article [Ol2] and the
references therein) underscores its significance. There is a remarkable
relation between the present theory of geometric structures and the moving
frame method. It turns out that the first is both a generalization and a
specialization of the second: It is a generalization because it incorporates
the concept of curvature which is not present in the second. It is a
specialization because when $\mathcal{R}=0,$ the first considers only
pseudogroups arising from transitive Lie group actions whereas the second
applies to much more general pseudogroups, not even transitive.

To clarify this relation, we fix some integer $1\leq k\leq n=\dim M,$ $p\in
M $ and define the fiber bundle $J_{k}(r,M)\rightarrow M$ whose fiber over $%
p\in M$ is the set consisting of the $k$-jets of (locally defined) maps $%
\mathbb{R}^{r}\rightarrow M$ with the source at the origin $o\in \mathbb{R}%
^{r},$ target at $p$ and have maximal rank $r$ at $o$ (hence near $o).\ $So
we have the obvious projections of fiber bundles

\begin{equation}
....\longrightarrow J_{k+1}(r,M)\longrightarrow J_{k}(r,M)\longrightarrow
....\longrightarrow J_{1}(r,M)\longrightarrow M
\end{equation}

We observe that the fiber of $J_{k}(r,M)\rightarrow M$ over $p$ can be
identified with the $r$-dimensional (local) submanifolds of $M$ passing
through $p$ modulo the equivalence relation defined by $k$-th order contact
at $p.$ Now the action of a transitive Lie group $G$ on $M$ lifts in the
obvious way to $J_{k}(r,M)$ for all $k$ and "stabilizes" at some order $m,$
i.e., $G$ acts reely on $J_{k}(r,M)$ for $k\geq m$ (For more general Lie
group and pseudogroup actions the stabilization is a more delicate problem,
see [AO], [OP]). This stabilization phenomenon which is first observed in
[Ov] and corrected and generalized in [Ol1], is a fundamental fact and is
the key to the theory of differential invariants by the moving frame method.

Once the action of $G$ becomes free on $J_{m}(r,M)$ with the orbit space $%
J_{m}(r,M)/\sim $ and the quotient map $\pi :J_{m}(r,M)\rightarrow
J_{m}(r,M)/\sim ,$ the moving frame method proceeds by choosing a crossection

\begin{equation}
c:J_{m}(r,M)/\sim \text{ \ }\longrightarrow J_{m}(r,M)
\end{equation}%
to the orbits of $G,$ i.e., $c$ chooses from each orbit a single element in
a smooth way. Since each orbit is in 1-1 correspondence with $G$ by freeness
and $J_{m}(r,M)$ is disjoint union of orbits, (109) gives a map

\begin{equation}
\widetilde{c}:J_{m}(r,M)\longrightarrow G
\end{equation}%
defined as follows: If $x\in J_{m}(r,M),$ then $x$ belongs to the orbit $\pi
(x)$ and $\widetilde{c}(x)$ is the unique element in $G$ which maps $x$ to $%
c(\pi (x)).$ It is easy to check that $\widetilde{c}$ commutes with the
action of $G$ (the choice of left/right moving frames arises if we work with
an abstract Lie group and dissappears for transformation groups). It is
crucial to observe that all orbits are equal in the moving frame method,
i.e., there is no canonical orbit and also there is no canonical crossection
(109). We also remark here that moving frame method is a \textit{local
theory and all spaces, maps...etc above have only local meaning with the
assumption that locally everything is "nice". }However, a global version is
proposed in [KL] for algebraic groups and their actions on algebraic
manifolds. Once the moving frame (110) is constructed, the moving frame
method allows us to compute everything in sight (and not in sight!)
algorithmically and constructively and we refer to the survey article [Ol2]
to give the reader an idea about the power and scope of this method.

Now we specialize to the case $r=\dim M=n.$ Fixing some $e\in M$
arbitrarily, the fiber of $J_{k}(n,M)\rightarrow M$ over $p$ can now be
identified with the fiber $\mathcal{U}_{k}^{e,p}$ of the principal bundle $%
\mathcal{U}_{k}^{e,\bullet }\rightarrow M$ over $p,$ i.e., $J_{k}(n,M)\cong 
\mathcal{U}_{k}^{e,\bullet }$ for all $k\geq 0.$ With this identification,
the prolonged action of $G$ on $J_{k}(n,M)$ becomes simply the above
mentioned composition with the arrows of $\mathcal{U}_{k}^{e,\bullet }$ at
the target, i.e., $g\in G$ maps $j_{k}(f)^{e,p}\in \mathcal{U}_{k}^{e,p}$ to 
$j_{k}(g)^{p,g(p)}\circ j_{k}(f)^{e,p}=j_{k}(g\circ f)^{p,g(p)}.$ Now
suppose we choose $m$ as in the first paragraph of the previous section.
Then $g$ fixes $j_{m}(f)^{e,p}$ $\iff j_{m}(f)^{e,p}=j_{m}(g)^{p,g(p)}\circ
j_{k}(f)^{e,p}\iff g(p)=p$ and $j_{m}(g)^{p,p}=Id\iff g=Id.$ Therefore $G$
acts freely on $\mathcal{U}_{m}^{e,\bullet }$ and the $m$ is the
stabilization number in the moving frame method. We observe that in this
special case $k=\dim M,$ there is also a canonical orbit which is $\mathcal{P%
}_{m}^{e,\bullet }\subset \mathcal{U}_{m}^{e,\bullet }$ where $\mathcal{P}%
_{m}\subset \mathcal{U}_{m}$ is defined by the action of $G$ on $M$ as in
the previous section. Therefore, we have the "moving frame" on this orbit

\begin{equation}
\mathcal{P}_{m}^{e,\bullet }\longrightarrow G
\end{equation}%
since any element of $g\in G$ is uniquely determined by its $m$-arrow from $%
e $ to $g(e).$ Observe that we already have the "restriction of the moving
frame $\mathcal{U}_{m}^{e,\bullet }\rightarrow G$ to the canonical orbit $%
\mathcal{P}_{m}^{e,\bullet }"$ whereas the moving frame $\mathcal{U}%
_{m}^{e,\bullet }\rightarrow G$ itself is not in sight as we have not
choosen any crossection yet !! Now let us take a closer look at the orbit
space $\mathcal{U}_{m}^{e,\bullet }/\sim .$ The stabilizer $H\subset G$ at $%
e $ acts on the fiber $\mathcal{U}_{m}^{e,e}\cong G_{m}(n)$ on the left and
the orbit space is the left coset space $\mathcal{U}_{m}^{e,e}/j_{m}(H).$
However, since $G$ acts transitively on $M,$ there is an obvious bijection

\begin{equation}
\mathcal{U}_{m}^{e,\bullet }/\sim \text{ \ \ }\iff \text{ \ \ \ }\mathcal{U}%
_{m}^{e,e}/j_{m}(H)\cong G_{m}(n)/j_{m}(H)
\end{equation}%
because any equivalence class in $\mathcal{U}_{m}^{e,\bullet }$ has a
representative in $\mathcal{U}_{m}^{e,e}$ and two arrows are related in $%
\mathcal{U}_{m}^{e,\bullet }$ if and only if "their projections" are related
in $\mathcal{U}_{m}^{e,e}.$ Therefore, the \textit{local }crossection (109)
amounts to choosing a local crossection

\begin{equation}
G_{m}(n)/j_{m}(H)\longrightarrow G_{m}(n)
\end{equation}%
around the coset defined by $j_{m}(H).$ In particular, if $\dim
G_{m}(n)/j_{m}(H)=t$ and $U\subset $ $G_{m}(n)/j_{m}(H)$ is a neighboorhood
of the coset $j_{m}(H)$ in $G_{m}(n)/j_{m}(H),$ then a coordinate
crossection amounts to introducing a coordinate patch

\begin{equation}
G_{m}(n)/j_{m}(H)\supseteqq U\longrightarrow \mathbb{R}^{t}
\end{equation}%
Now we compare (107) and (114). There are two main differences: 1) The order
is $m+1$ in (107) because our approach to geometric structures incorporates
splittings and therefore the concept of curvature. 2) Since $%
G_{m}(n)/j_{m}(H)$ is a smooth manifold, we can always seperate the orbits
locally by smooth functions as in (114). However, to seperate the orbits
globally by polynomial functions as in (107), we need an affirmative answer
to \textbf{Q1, Q2.}

A last remark: As we have seen, our approch to geometric structures arises
from the special case of (108) for $r=\dim M.$ The key property of the tower
(108) is that any transitive Lie group action on $M$ prolongs to this tower
and stabilizes at some order. There are many other such (nonlinear and
linear) towers. Calling such a tower a jet-representation, the theory of
PHG's can be developed in the framework of jet-representations, the PHG's
themselves arising from the above canonical jet-representation with $r=\dim
M $ ([Or3]).

\bigskip

\textbf{References}

\bigskip

[AO] Adams, S., Olver, P., Prolonged analytic connected group actions are
generalically free, Transformation Groups, 23, (2018), 893-913

[Bl] Blaom, A., Geometric structures as deformed infinitesimal symmetries,
Trans. Amer. Math Soc. 358, 2006, 3651-71

[CS] Cap,A., Slovak,J., Parabolic Geometries I : Background and General
Theory, Mathematical Surveys and Monographs, 154, American Mathematical
Society, 2009

[D] Draisma,J., Lie Algebras of Vector Fields, PhD Thesis, Technische
Universiteit, Eindhoven, Netherlands, 2002

[FO1] Mark Fels, Peter J. Olver: Moving Coframes. I. A practical algorithm,
Acta Appl. Math. 51 (1998), 161-213

[FO2] Mark Fels, Peter J. Olver: Moving Coframes. II. Regularization and
theoretical foundations, Acta Appl. Math. 55 (1999), 127-208

[KL] Kruglikov,B., Lychagin,V., Global Lie-Tresse Theorem, Sel. Math. New
Ser. 22 (2016), 1357-1411

[Or1] Orta\c{c}gil, E.H., An Alternative Approach to Lie Groups and
Geometric Structures, Oxford University Press, 2018

[Or2] Orta\c{c}gil, E., Riemannian geometry as a curved prehomogeneous
geometry, arXiv: 1003.3220, 2010

[Or3] Orta\c{c}gil, E., Representations of prehomogeneous geometries, in
progress

[Ol1] Olver, P.J., Equivalence, Invariants and Symmetry, Cambridge
University Press, 1995

[Ol2] Olver, P.J., Lectures on moving frames, in: Symmetries and
Integrability of Difference Equations, D. Levi, P. Olver, Z. Thomova, and P.
Winternitz, eds., London Math. Soc. Lecture Note Series, vol. 381, Cambridge
University Press, Cambridge, 2011, 207-246

[OP] Olver, P.J., and Pohjanpelto, J., Persistence of freeness for
pseudo-group actions, Arkiv Mat. 50 (2012), 165-182.

[Ov] Ovsiannikov, L.V., Group Analysis of Differential Equations, Academic
Press, New York, 1982

[Po] Pommaret, J.F., Partial Differential Equations and Group Theory, New
Perspectives for Applications, Kluwer Academic Publishers, 1994

\bigskip

Erc\"{u}ment Orta\c{c}gil

Bodrum, ortacgile@gmail.com

\end{document}